\documentclass[11pt,a4paper,reqno]{amsart}

\usepackage{times} % assumes new font selection scheme installed
\usepackage{amsmath} % assumes amsmath package installed
\usepackage{amssymb}  % assumes amsmath package installed
\usepackage{amsthm}
\usepackage{latexsym}
\usepackage{amsfonts}
\usepackage{xcolor}
\usepackage{mathtools}
\usepackage{enumerate}
\usepackage{cite}
\usepackage{graphicx}

\newtheorem{thm}{Theorem}[section]

\newtheorem{lem}[thm]{Lemma}
\newtheorem{prop}[thm]{Proposition}
\theoremstyle{definition}
\newtheorem{defn}[thm]{Definition}
\newtheorem{rem}[thm]{Remark}
\newtheorem{ex}[thm]{Example}

\numberwithin{equation}{section}
\allowdisplaybreaks
\newcommand{\setdef}[2]{\left\{ #1 \left\vert\vphantom{#1} #2 \right.\right\}}

\newcommand{\Id}{I}%\operatorname{Id}}

\renewcommand{\Re}{\operatorname{Re}}

\newcommand{\ran}{\operatorname{ran}}
\newcommand{\dom}{\operatorname{dom}}

\newcommand{\wt}{\widetilde}

\newcommand{\<}{\langle}
\renewcommand{\>}{\rangle}

\newcommand{\R}{\ensuremath{\mathbb R}}    % Reelle Zahlen
\newcommand{\C}{\ensuremath{\mathbb C}}    % Komplexe Zahlen
\newcommand{\N}{\ensuremath{\mathbb N}}    % Nat"urliche Zahlen
         \newcommand{\frakA}{\mathfrak A}
         \newcommand{\frakB}{\mathfrak B}
         \newcommand{\frakC}{\mathfrak C}
         \newcommand{\frakD}{\mathfrak D}
         
\newcommand{\calF}{\mathcal F}         
\newcommand{\calG}{\mathcal G}         
\newcommand{\calH}{\mathcal H}

\newcommand{\calR}{\mathcal R}

\newcommand{\calU}{\mathcal U}         
\newcommand{\calV}{\mathcal V}         
\newcommand{\calW}{\mathcal W}         
\newcommand{\calX}{\mathcal X}         
\newcommand{\calY}{\mathcal Y}         
\newcommand{\calZ}{\mathcal Z}         

\newcommand{\la}{\lambda}

\newcommand{\bmat}[4]
{
   \begin{bmatrix}
      #1 & #2\\
      #3 & #4
   \end{bmatrix}
}
\newcommand{\bvek}[2]
{
   \begin{bmatrix}
      #1\\
      #2
   \end{bmatrix}
}

\newcommand{\sbvek}[2]{\left[\begin{smallmatrix}#1\\#2\end{smallmatrix}\right]}
\newcommand{\spvek}[2]{\left(\begin{smallmatrix}#1\\#2\end{smallmatrix}\right)}
\newcommand{\sbmat}[4]{\left[\begin{smallmatrix}#1 & #2\\#3 & #4\end{smallmatrix}\right]}

\renewcommand{\div}{\operatorname{div}}
\DeclareMathOperator{\curl}{curl}

\title[Infinite-dimensional port-Hamiltonian systems]{Infinite-dimensional port-Hamiltonian systems  - a system node approach}

\author[F.\ Philipp]{Friedrich Philipp}
\address{{\bf F.~Philipp:} Institute for Mathematics, Faculty of Mathematics and Natural Sciences, Technische Universit\"at Ilmenau, Ilmenau, Germany}
\email{friedrich.philipp@tu-ilmenau.de}
\urladdr{https://www.tu-ilmenau.de/obc/team/friedrich-philipp}

\author[T.\ Reis]{Timo Reis}
\address{{\bf T.~Reis:} Institute for Mathematics, Faculty of Mathematics and Natural Sciences, Technische Universit\"at Ilmenau, Ilmenau, Germany}
\email{timo.reis@tu-ilmenau.de}
\urladdr{https://www.tu-ilmenau.de/systpde/team/timo-reis}

\author[M.\ Schaller]{Manuel Schaller}
\address{{\bf M.~Schaller:} Institute for Mathematics, Faculty of Mathematics and Natural Sciences, Technische Universit\"at Ilmenau, Ilmenau, Germany}
\email{manuel.schaller@tu-ilmenau.de}
\urladdr{https://www.tu-ilmenau.de/obc/team/manuel-schaller}

\begin{document}
\begin{abstract}
We consider an operator-theoretic approach to linear infinite-dimensional port-Hamiltonian systems. In particular, we use the theory of system nodes
by {\sc Staffans} \cite{Staffans2005} to formulate a~suitable concept for port-Hamiltonian systems, which allows a unifying approach to systems with
boundary as well as distributed control and observation. The concept presented in this article is further neither limited to parabolic nor hyperbolic
systems, and it also covers partial differential equations on multi-dimensional spatial domains. Our presented theory is substantiated by means of several
physical examples.

\end{abstract}

\maketitle

\section{Introduction}

Port-Hamiltonian systems arise from energy-based modeling of physical systems.
This class is further closed under power-conserving network interconnection, that is, loosely speaking, the coupling of port-Hamiltonian systems again
results into a port-Hamiltonian system. Hence, they constitute a~network modeling paradigm that
offers a systematic approach for the modeling of interactions of systems from various physical domains, such as electrodynamics, (continuum) mechanics, flow
dynamics, chemical reaction kinetics, and thermodynamics.\\
Theory has made tremendous progress for port-Hamiltonian systems governed by ordinary dif\-ferential and differen\-tial-algebraic equations
\cite{beattie2018linear,van2014port,gernandt2021linear,MS22a,MS22b}. On the other hand, many partial differential equations arising in modelling of physical
systems have intrinsically a structure which apparently resembles the structure that is well-known from the case of ordinary differential equations. As
partial differential equations can usually be rewritten as infinite-dimensional, a suitable operator-theoretic approach is needed for an analytically sincere
proceeding to this kind of port-Hamiltonian systems, which should preferably neither be limited to distributed nor to boundary control, and also not be
restricted to special classes of partial differential equations, such as, for instance, parabolic or hyperbolic ones.

The most general approach  is via geometric structures such as  Dirac structures and Lagrange submanifolds, which are essential for the availability of an
energy balance. The aforementioned geometric structures are, roughly speaking, related to bundles of vector spaces which are self-orthogonal with respect to certain
indefinite inner products. This altogether gives rise to a~behavioral approach to port-Hamiltonian systems, see \cite{van2014port} for an introduction.
This has been generalized to classes of partial differential equations with boundary ports in
\cite{van2021differential,van2002hamiltonian,schoberl2014jet,le2005dirac}, where the concept of Stokes-Dirac structures has been used, which is a~special
Dirac structure on function spaces. A~functional analytic approach to Dirac structures and Lagrange submanifolds is presented
in \cite{reis2021some}. The disadvantage of the geometric approaches to port-Hamiltonian systems is that they provide neither existence results for solutions
nor qualitative characterizations of the solutions, such as stability or regularity.

For spatially one-dimensional hyperbolic partial differential equations such a rigorous analysis is presented in initially in the thesis \cite{Vill07}, and
further elaborated in the textbook \cite{Jacob2012} and several subsequent publications \cite{JMZ19,JK19,JKZ21}, where infinite-dimensional systems theory
has been successfully applied to carry out a deep qualitative analysis of such systems. In particular, the one-dimensionality of the spatial domain leads to the
consideration of a finite-dimensional boundary trace space, which allows the formulation of handy linear algebraic criteria on the boundary conditions for
existence and stability of solutions. Quite recently, a~very promising generalization to hyperbolic partial differential equations on multidimensional
spatial domains has been introduced \cite{Skrepek21}.

The approach in this article will be somewhat more operator theoretic, but mainly motivated by the representation
\begin{equation}\label{e:phs}
\begin{aligned}
\dot x(t)&=(J\;\:-\;R\:)Hx(t)+(B-P)u(t),\\
y(t)&=(B^*+P^*)Hx(t)+(S-N)u(t)
\end{aligned}
\end{equation}
of finite-dimensional port-Hamiltonian systems from \cite{beattie2018linear,MehrMora19}, where
$J\in\C^{n\times n}$ and $N\in\C^{m\times m}$ are skew-Hermitian, and $H\in\C^{n\times n}$ and $W := \sbmat RP{P^*}S\in\C^{(n+m)\times (n+m)}$ are Hermitian
positive semi-definite.
The total energy of the system is given by the {\em Hamiltonian} $\calH(x) = \frac 12\, x^* H {x}$. By using \eqref{e:phs} together with skew-adjointness
of $J$ and $N$, as well as positive semi-definiteness of $W$ and $H$, for all $t>0$, $u\in L^2([0,t];\C^m)$, $x_0\in\C^n$, the solution of
\eqref{e:phs} with $x(0)=x_0$ fulfills the {\em dissipation inequality}
\begin{multline}
\calH(x(t))-\calH(x_0)=-\Re\int_0^t \spvek {Hx(\tau)}{u(\tau)}^*\sbmat RP{P^*}S\spvek {Hx(\tau)}{u(\tau)}{\rm d}\tau+\Re\int_0^t  u(\tau)^*{y(\tau)}{\rm d}\tau\\
\leq \Re\int_0^t u(\tau)^* {y(\tau)}{\rm d}\tau,\label{eq:enbalfin}
\end{multline}
which has the physical interpretation of an energy balance. Namely, whilst $\calH(x(t))$ stands for the energy stored at time $t$, the first integral after the equality sign is the energy dissipated by the system during the time interval $[0,t]$, whereas $\Re(u(\tau)^*{y(\tau)})$ can be regarded as the external power supply to the system.

Though it is claimed in \cite{MehrMora19} that the definition for the finite-dimensional case can be extended to state spaces with infinite
dimension,
the situation in infinite dimensions is by far not that simple, in particular when aiming to include
physically important partial differential equations with boundary control and observation. Hence, a~direct generalization of \eqref{e:phs} to the
infinite-dimensional case by replacing all matrices by (possibly unbounded) operators (which indeed would not be mathematically challenging) leads to a~class that leaves out a~variety of physical systems though these are port-Hamiltonian in some sense.\\
Our~concept for port-Hamiltonian systems is based on rewriting \eqref{e:phs} into
\begin{align}\label{e:phs_compact}
\spvek{\dot x(t)}{y(t)} &= \sbmat{\Id_n}{}{}{-\Id_m} M\spvek{Hx(t)}{u(t)},\\[1mm]
 M&:=\sbmat{J-R}{B-P}{-B^*-P^*}{N-S},\label{e:phs_compact2}
\end{align}
where, by the assumptions on $J$, $R$, $B$, $P$, $N$ and $S$, the matrix $ M\in\C^{(n+m)\times (n+m)}$ is dissipative in the sense that
$ M+ M^*$ is negative semi-definite.\\
Our concept covering the infinite-dimensional case will indeed be based on
a surrogate for the representation~\eqref{e:phs_compact} involving some dissipative operator $ M$ and some~nonnegative operator $H$. To include partial differential
equations with boundary control and observation, $ M$ may be unbounded and defined on some proper subspace which is not necessarily a~Cartesian
product of the state and input %and output
space. Hereby, we make use of the rich theory of {\em system nodes}  by {\sc Staffans}
 \cite{Staffans2005}, that offers a wide flexibility in modeling boundary control systems. This framework is a
generalization of the class of well-posed systems in the sense of \cite{Staffans2005}. One feature of system nodes is that the control and observation
equations take the form
\begin{equation}\label{eq:sysnode}
\spvek{\dot{x}(t)}{y(t)}=\sbvek{A\&B\\[-1mm]}{C\&D}\spvek{x(t)}{u(t)}
\end{equation}
rather than a split $\dot{x}(t)=Ax(t)+B{u(t)}$ and
$y(t)=Cx(t)+D{u(t)}$. Hereby, $A\&B$ and $C\&D$ may be a composite operators defined on some dense subspace of the Cartesian product of the state and input
space. This allows for a natural, simple and direct modeling of boundary control.
By a~suitable combination of \eqref{e:phs} with \eqref{eq:sysnode}, we will present a novel and analytically sound approach to linear infinite-dimensional port-Hamiltonian systems.

If $H$ is the identity, we obtain a~class which is called {\em internally impedance passive} in \cite{staffans2002passive}, and this system type
has also been briefly mentioned in \cite{Vill07} in the context of port-Hamiltonian systems. In this case, the Hamiltonian is simply given by
$\calH(x)=\tfrac12\|x\|^2$. The incorporation of nontrivial Hamiltonians, which do not necessarily have to be coercive or bounded, will be shown to enable a
by far wider class of energy functionals.

This article is organized as follows: The mathematical background is briefly presented in Section~\ref{sec:prelim}. We introduce notation and the basics on
system nodes and (maximal) dissipative operators. Thereafter, Section~\ref{sec:phsysnodes} contains the main part of this article: We present the system node approach to port-Hamiltonian systems. The operator $H$ will be constructed by means of closed symmetric sesquilinear forms, whereas $ M$ will be introduced as a~dissipative operator with certain properties. Some further results on solvability will be shown. In Section~\ref{sec:appl}, we will consider a~variety of practical examples that fit into our framework. %Finally in Section~\ref{sec:passive}, we will present some relations between port-Hamiltonian systems and systems which are impedance passive in the sense of \cite{staffans2002passive}.

\section{Mathematical preliminaries}\label{sec:prelim}

\subsection{Notation}

Let $\calX$, $\calY$ be Hilbert spaces,
which are always assumed to be complex throughout this work.
The norm in {$\calX$} will be denoted by $\|\cdot\|_{{\calX}}$ or {simply} $\|\cdot\|$, if clear from context. The space $\C^n$ will be always equipped with the standard Euclidean inner product. The identity mapping in $\calX$ is abbreviated by $\Id_{\calX}$ (or
just $\Id$, if clear from context), and we set $\Id_n:=\Id_{\C^n}$.

The symbol $\calX^*$ stands for the {\em anti-dual} of $\calX$, that is, the space of all continuous antilinear (i.e., additive and conjugate homogeneous) functionals on $\calX$.
 Hence, the canonical duality product $\<\cdot,\cdot\>_{\calX^*,\calX}$ is (as well as the inner product $\<\cdot,\cdot\>_{\calX}$ in $\calX$) a~sesquilinear form. That is, it is linear in the first argument, and antilinear in the second argument. Further, note that
 the Riesz map $\calR_\calX$, sending $x\in\calX$ to the functional $\<x,\cdot\>_{\calX}$ is a linear isometric isomorphism from $\calX$ to $\calX^*$. Again, if the spaces are clear from context, we may skip the subindices. Further, if not stated otherwise, a~Hilbert space is canonically identified with its anti-dual. Note that, in this case, $\calR_\calX=\Id_\calX$.

 The space of bounded linear operators from $\calX$ to $\calY$ is denoted by $L(\calX,\calY)$. As usual, we abbreviate $L(\calX):= L(\calX,\calX)$. The domain $\dom A$ of a possibly unbounded linear operator $A:\calX\supset\dom A\to\calY$ is typically equipped with the graph norm $\|x\|_{\dom A}:=\big(\|x\|_{\calX}^2+\|Ax\|_{\calY}^2\big)^{1/2}$. %The closure of a~closable operator $A:\calX\supset\dom A\to\calY$ is denoted by $\overline{A}$, and for two operators $A:\calX\supset\dom A\to\calY$, $B:\calX\supset\dom B\to\calY$ we write $A\subset B$, if $\dom A\subset\dom B$, and $Ax=Bx$ for all $x\in\dom A$. In this case, we call $B$ an \textcolor{blue}{``extension of $A$''}. Equivalently, we call $A$ a~``restriction of $B$''.

 The adjoint $A^*:\calY^*\supset\dom A^*\to\calX^*$ of a~densely defined linear operator $A:\calX\supset\dom A\to\calY$ has the domain
\[
\dom A^*=\setdef{y'\in \calY^*}{\exists\, z'\in\calX^*\text{ s.t.\ }\forall\,x\in\dom A:\;\langle y',Ax\rangle_\calY=\langle z',x\rangle_\calX }.
\]
The functional $z'\in \calX^*$ in the above set is uniquely determined by $y'\in\dom A^*$, and we set $A^*y'=z'$. Note that we identify $\C^{n\times m}\cong L(\C^m,\C^n)$. Together with the fact that $\C^n$ and $\C^m$ are equipped with the Euclidean inner product, this means that $A^*\in\C^{n\times m}$ is the conjugate transpose of $A\in\C^{m\times n}$. Likewise, $x^*$ is the conjugate transpose of $x\in\C^n\cong \C^{n\times 1}$, such that the inner product in $\C^n$ reads
\[
\langle x,y\rangle_{\C^n}=y^*x.
\]
 The resolvent set of $A: \calX\supset\dom A\to \calX$ is denoted by $\rho(A)$, i.e., \[\rho(A) = \setdef{\la\in\C}{(\la \Id - A)^{-1}\in L(\calX)}.\]
% The identity operator on the space $\calX$ is $\Id_\calX$ %The symbol $A^*$ stands for the adjoint of a~linear operator $A$.
We denote the open right half-plane by $\C_+ = \setdef{\la\in\C}{\Re\la > 0}$ and the set of the nonnegative real numbers by $\R_{\ge0}$.% We will frequently denote the Cartesian product of two sets $A$ and $B$ by $\sbvek AB$.

We use the notation of the book \cite{adams2003sobolev} by {\sc Adams} for Lebesgue and Sobolev spaces as well as for the space of $k$ times continuously differentiable functions. For function spaces with values in a~Hilbert space $\calX$, we indicate this by denoting ``$;\calX$'' after specifying the (spatial or temporal) domain. For instance, the Lebesgue space of $p$-integrable $\calX$-valued functions on the domain $\Omega$ is $L^p(\Omega;\calX)$.
Note that, throughout this article, integration of $\calX$-valued functions always has to be understood in the Bochner sense \cite{Diestel77}.

%  of $k$-times differentiable $\calX$-valued functions on an interval $J$  will be denoted by $C^k(J;\calX)$, and the space of continuous functions by $C^0(J;\calX)=C(J;\calX)$. Further, for $p\in [1,\infty]$, $k\in\N$, the Lebesgue space $L^p$ and Sobolev space $W^{k,p}$ of $\calX$-valued functions on $J$ are denoted by
%$W^{k,p}(J;\calX)$ and $L^p(J;\calX)$. %For a~subinterval $J_0\subset J$ we canonically identify $C^k(J_0;\calX)\subset C^k(I;\calX)$, $W^{k,p}(J;\calX)\subset W^{k,p}(I;\calX)$ and $L^p(J;\calX)\subset L^p(I;\calX)$ via restriction of functions to $J$.
% Further, $W^{k,p}_{\rm loc}(J;\calX)$ and $L^p_{\rm loc}(J;\calX)$ are the spaces of functions $f:I\to \calX$ with $f\in L^p(K;\calX)$ ($f\in W^{k,p}(K;\calX)$) for all compact subintervals $K\subset I$.

%\marginpar{should we add direct sum $\oplus$?}

\subsection{Operator and system nodes}
Let $\calX$, $\calU$, and $\calY$ be Hilbert spaces and denote the canonical projection onto $\calX$ ($\calY$) in $\calX\times\calY$ by $P_\calX$ ($P_\calY$, respectively). Let
\[S : \calX\times\calU\supset\dom S\to\calX\times\calY\]
be a linear operator. Its so-called {\em main operator} $A : \calX\supset\dom A\to\calX$ maps from $\dom A := \setdef{x\in\calX}{\spvek x0\in\dom S}$ to $\calX$ via $Ax := P_\calX S\spvek x0$. {We set}
\[A\& B := P_\calX S
\qquad\text{and}\qquad
C\& D := P_\calY S,
\]
so that $S$ can be written as
\[
S = \sbvek{A\& B}{C\& D}.
\]
%We first provide the modeling framework for the control systems that we are going to consider in this paper. For this, let $\calX$, $\calU$, and $\calY$ be Hilbert spaces and consider a linear operator $S:\calX\times\calU\supset\dom S\to \calX \times \calY$ describing the (abstract) dynamics
The concept of operator nodes poses natural assumptions on the operator $S$, in order to guarantee favorable properties and a suitable solution concept to the abstract dynamics
\begin{equation}
\spvek{\dot{x}(t)}{y(t)}
= S \spvek{{x}(t)}{u(t)}.\label{eq:ODEnode}\end{equation}
We provide the definition from \cite[Def.~2.1]{OpmeerStaffans2014}. For basics on strongly continuous semigroups on Hilbert spaces we refer to \cite{TuWe09}.

\begin{defn}[Operator/system node]\label{def:opnode}
An {\em operator node} on the triple $(\calY,\calX,\calU)$ is a (possibly unbounded) linear operator $S : \calX\times\calU\supset\dom S\to\calX\times\calY$ with the following properties:
\begin{enumerate}[(a)]
\item\label{def:opnode1} $S$ is closed;
\item\label{def:opnode2} $P_\calX S:\calX\times\calU\supset\dom S\to\calX$ is closed;
\item\label{def:opnode3} for all $u\in \calU$, there exists some $x\in \calX$ with $\spvek{x}{u}\in \dom S$;
\item\label{def:opnode4} its main operator $A$ is densely defined and has a nonempty resolvent set.
\end{enumerate}
If, further, $A$ generates a strongly continuous semigroup on $\calX$, then $S$ is called {\em system node}.
\end{defn}

\begin{rem}[Operator/system nodes I]\label{rem:nodes}
Let $S = \sbvek{A\& B}{C\& D}$ be an~operator node on $(\calY,\calX,\calU)$.
\begin{enumerate}[(a)]
%\item\label{rem:nodes1} It follows from the definition of operator nodes that $C\&D\in \mathcal L(\dom(A\&B),\calY)$, where $\dom(A\&B)$ is equipped with the graph norm of $A\&B$. In particular, the operator $C$ with $Cx:= C\&D\spvek x0$ fulfills $C\in \mathcal L(\dom A,\calY)$.
\item\label{rem:nodesa} By $\calX_{-1}$ denote the completion of $\calX$ with respect to the norm $\|x\|_{\calX_{-1}}:= \|(\alpha \Id-A)^{-1}x\|$ for some $\alpha\in\rho(A)$. Note that the topology of $\calX_{-1}$ does not depend on the particular choice of $\alpha\in \rho(A)$ \cite[Prop.~2.10.2]{TuWe09}. Then the operator $A$ extends to closed and densely defined operator $A_{-1} : \calX_{-1}\supset\dom A_{-1} = \calX\to\calX_{-1}$ such that $A_{-1} : \calX\to\calX_{-1}$ is bounded. The spectra of $A$ and $A_{-1}$ coincide. Moreover, if $S$ is a system node, the semigroup $\frakA(\cdot)$ generated by $A$ extends to a~semigroup $\frakA_{-1}(\cdot)$ on $\calX_{-1}$. The generator of this semigroup is $A_{-1}$ \cite[Prop.~2.10.3 \& 2.10.4]{TuWe09}.

\item\label{rem:nodesb} There exists an operator $B\in L(\calU,\calX_{-1})$ such that $[A_{-1}\ B]\in L(\calX\times\calU,\calX_{-1})$ is an extension of $A\& B$. The domain of $A\&B$ (equally: the domain of $S$) satisfies
\[
\dom(A\&B)=\setdef{\spvek xu \in \calX\times \calU}{A_{-1}x+Bu\in \calX},
\]
see \cite[Def.~4.7.2 \& Lem.~4.7.3]{Staffans2005}.

\item\label{rem:nodesc} Given $u\in\calU$,
%By using that generators of semigroups are densely defined (see \cite[Chap.~2, Thm.~1.5]{engel2000one}),
the affine space
\[
\setdef{x\in \calX}{\spvek xu\in \dom(A\& B)}
\]
is dense in $\calX$.

\item\label{rem:nodese} For all $\alpha\in\rho(A)$, the norm
\[\left\|\spvek xu\right\|_\alpha:= \left(\|x-(\alpha \Id-A_{-1})^{-1}Bu\|_\calX^2+\|u\|_\calU^2\right)^{1/2}\]
is equivalent to the graph norm of $S$. Moreover, the operator
\[\sbmat \Id{-(\alpha \Id -A_{-1})^{-1}B}0\Id\]
maps $\dom S$ bijectively to $\dom A\times \calU$ \cite[Lem.~4.7.3]{Staffans2005}.
\end{enumerate}
\end{rem}

%It follows from (iii) and (iv) that an operator node $S$ is always densely defined. Hence, $S^*$ exists as a closed and densely defined operator and it can be shown (see \cite[Proposition 2.4]{malinen2006linear}) that $S^*$ is an operator node on $(X,Y,U)$ with main operator $A^*$.
Now we record an~operator theoretic lemma which is needed at several places. Analogous to the statements at the beginning of this section, for Hilbert spaces $\calV$, $\calW$, $\calZ$, and an operator $R:\calZ\supset\dom R\to \calV\times\calW$,
we denote the canonical projection onto $\calV$ ($\calW$) in $\calV\times\calW$ by $P_\calV$ ($P_\calW$, respectively).

\begin{lem}[Closedness of block operators]\label{lem:closedblock}
Let $\calV$, $\calW$, $\calZ$ be Hilbert spaces, and let $R:\calZ\supset\dom R\to \calV\times\calW$ be a~linear operator. Then the following statements are equivalent:
\begin{enumerate}[(i)]
\item\label{lem:closedblock1} $R$ is closed, and $P_\calV R$ is closed with domain $\dom(P_\calV R)=\dom R$.
\item\label{lem:closedblock2} $P_\calV R$ is closed with domain $\dom(P_\calV R)=\dom R$, and $P_\calW R\in L(\dom(P_\calV R),\calW)$.
\item\label{lem:closedblock3} $R$ is closed, and $P_\calW R\in L(\dom(P_\calV R),\calW)$.
\end{enumerate}
\end{lem}
\begin{proof}\

\noindent``\eqref{lem:closedblock1}$\Rightarrow$\eqref{lem:closedblock2}'': Assume that \eqref{lem:closedblock1} is true. We show $P_\calW R\in L(\dom(P_\calV R),\calW)$ by employing the closed graph theorem \cite[Thm.~7.9]{Alt16}.
Assume that
\[
(z_n)\to z\text{ in $\dom(P_\calV R)$},\qquad
(P_\calW R z_n)\to w\text{ in $\calW$.}
\]
%By closedness of {$P_\calV R$}, we have that $z\in\dom(P_\calV R)=\dom R$.
This implies that $(R z_n)$ converges to $\spvek{P_\calV R z}{w}\in\calV\times\calW$, and we obtain from the closedness of $R$ that $R z=\spvek{P_\calV R z}{w}$. In particular $P_\calW Rz=w$.\\
``(ii)$\Rightarrow$(iii)'': Under assumption (ii), consider a~sequence with
\[
(z_n)\to z\text{ in $\calZ$},\qquad
(R z_n)\to \spvek{v}w\text{ in $\calV\times \calW$.}
\]
Then $(P_\calV R z_n)$ converges in $\calV$ to $v$, and closedness of $P_\calV R$ with domain $ \dom R$ leads to $z\in\dom R$ with $P_\calV R z=v$. The property $P_\calW R\in L(\dom (P_\calV R),\calW)$ further leads to convergence of $(P_\calW R z_n)$ to $P_\calW R z$. Altogether, $z\in\dom R$ with $Rz =\spvek{v}w$.\\
``\eqref{lem:closedblock3}$\Rightarrow$\eqref{lem:closedblock1}'': Assume that $R$ is closed, and $P_\calW R\in L(\dom(P_\calV R),\calW)$, and let
\[
(z_n)\to z\text{ in $\calZ$},\qquad
(P_\calV R z_n)\to v\text{ in $\calV$.}
\]
Then $(P_\calW Rz_n)$ is a~Cauchy sequence in $\calW$, and thus convergent to some $w\in\calW$. Hence, $(R z_n)$ converges to $\spvek{v}w$, and closedness of $R$ leads to $z\in\dom R$ and $R z=\spvek{v}w$. In particular, $P_\calV R z=w$.
\end{proof}

\begin{rem}[Operator/system nodes II]\label{rem:nodesnew}
Let $S = \sbvek{A\& B}{C\& D}$ be an~operator node on $(\calY,\calX,\calU)$.
Lemma~\ref{lem:closedblock} yields that $C\&D\in \mathcal L(\dom(A\&B),\calY)$, where $\dom(A\&B)$ is equipped with the graph norm of $A\&B$. In particular, for the operator $C$ with $Cx:= C\&D\spvek x0$, we have $C\in \mathcal L(\dom A,\calY)$.\\
On the other hand, if $S = \sbvek{A\& B}{C\& D}$ has properties \eqref{def:opnode2}-\eqref{def:opnode3} in Definition~\ref{def:opnode} and, additionally, $C\&D\in \mathcal L(\dom(A\&B),\calY)$, then Lemma~\ref{lem:closedblock} yields that $S$ is an~operator node.
\end{rem}
The previous two remarks allow to define the concept of the transfer function.

\begin{defn}[Transfer function]
Let $S = \sbvek{A\& B}{C\& D}$ be a~system node on $(\calY,\calX,\calU)$. The \emph{transfer function $\calG$ associated with $S$} is
\[
\begin{aligned}
\calG\colon &&\rho(A)\to &\,\mathcal L(\calU,\calY),\\&&s\mapsto&\,C\& D \sbvek{(s\Id-A_{-1})^{-1}B}\Id.
\end{aligned}
\]
\end{defn}

\begin{rem}[Adjoint operator node]\label{rem:adjnode}
Let $S = \sbvek{A\& B}{C\& D}$ be an~operator node on $(\calY,\calX,\calU)$.
It has been shown in \cite[Prop.~2.4]{malinen2006linear} that the adjoint $S^*:\calX\times\calY\supset\dom S^*\to\calX\times\calU$ has the form
    \[S^*=\left[\begin{smallmatrix}{[A\& B]^d}\\{[C\& D]^d}\end{smallmatrix}\right]\]
with domain
\[\dom S^*=\setdef{\spvek xy\in \calX\times\calY}{A^*_{-1}x+C^*x\in \calX},\]
where $[A\&B]^d$ is the restriction of $[(A^*)_{-1}\; C^*]$ to $\dom S^*$, and, for all $\alpha\in\rho(A)$,
\[
[C\&D]^d\spvek xy:=B^*\Big(x-(\overline{\alpha}\Id-(A^*)_{-1}\big)^{-1}C^*\Big)
+\calG(\alpha)^*y.
\]
In particular, the main operator of $S^*$ is $A^*$.\\
If, additionally, $A$ generates a~strongly continuous semigroup $\frakA(\cdot)$ on $\calX$, then $A^*$ generates the adjoint semigroup $\frakA^*(\cdot)$ \cite[Prop.~2.8.5]{TuWe09}. Consequently, if $S$ is even a~system node, then $S^*$ is a~system node as well.
\end{rem}

Next we briefly recall suitable solution concepts for the differential equation
\eqref{eq:ODEnode} with $S=\sbvek{A\&B}{C\&D}$ being a system node.

\begin{defn}[Classical/generalized trajectories]\label{def:traj}
Let $T>0$, and let $S = \sbvek{A\& B}{C\& D}$ be a~system node  on $(\calY,\calX,\calU)$.\\
A {\em classical trajectory} for \eqref{eq:ODEnode} on $[0,T]$ is a triple
\[
(y,x,u)\,\in\, C([0,T];\calY)\times C^1([0,T];\calX)\times C([0,T];\calU)
\]
which for all $t\in[0,T]$ satisfies \eqref{eq:ODEnode}.\\
A {\em generalized trajectory} for \eqref{eq:ODEnode} on $[0,T]$ is a triple
\[
(y,x,u)\,\in\,L^2([0,T];\calY)\times C([0,T];\calX)\times  L^2([0,T];\calU),
\]
which is a~limit of classical trajectories for \eqref{eq:ODEnode} on $[0,T]$ in the topology of $L^2([0,T];\calY)\times C([0,T];\calX)\times L^2([0,T];\calU)$.
\end{defn}

If $S = \sbvek{A\& B}{C\& D}$ is a~system node  on $(\calY,\calX,\calU)$, then $A\&B$ can be regarded as a~system node on $(\{0\},\calX,\calU)$. Consequently, we may further speak of classical and generalized trajectories $(x,u)$ for
\begin{align}\label{e:ODE}
\dot{x}=A\&B\spvek xu.
\end{align}

The following result ensures the existence of unique classical trajectories with suitable control functions and initial values.

\begin{prop}[Existence of classical trajectories {\cite[Thm.~4.3.9]{Staffans2005}}]\label{prop:solex}
Let $S$ be a system node on $(\calY,\calX,\calU)$, let $T{>0}$, $x_0\in \calX$ and $u\in W^{2,1}([0,T];\calU)$ with $\spvek{x_0}{u(0)}\in \dom S$. Then there exists a unique
classical trajectory $(y,x,u)$ for \eqref{eq:ODEnode} with $x(0)=x_0$.
\end{prop}

We provide some further statements on classical generalized solutions.

\begin{rem}[Classical/generalized trajectories]\label{rem:sols}
Let $T{>0}$, and let $S = \sbvek{A\& B}{C\& D}$ be a~system node  on $(\calY,\calX,\calU)$.
\begin{enumerate}[(a)]
\item\label{rem:sols1} Assume that $(x,u)$ is a~classical trajectory for \eqref{e:ODE}. Then
\[
\spvek xu\in C([0,T];\dom S).
\]
\item\label{rem:sols2} $(x,u)$ is a~generalized trajectory for \eqref{e:ODE} if, and only if, $x\in C([0,T];\calX)$ and
\begin{equation}\label{eq:mildsol}
\forall\,t\in[0,T]:\quad x(t)=\frakA(t)x(0)+\int_0^t \frakA_{-1}(t-\tau)Bu(\tau){\rm d}\tau,
\end{equation}
where the latter has to be interpreted as an integral in the space $\calX_{-1}$ with $B\in L(\calU,\calX_{-1})$ as in Remark \ref{rem:nodes} (\ref{rem:nodesb}). %In particular, $x\in C([0,T];\calX_{-1})$.

\item\label{rem:sols3} If $(x,u,y)$ is a~generalized trajectory for \eqref{eq:ODEnode}, then, clearly, $(x,u)$ is a~generalized trajectory for \eqref{e:ODE}. In particular, \eqref{eq:mildsol} holds. The output evaluation $y(t)=C\&D \spvek {x(t)}{u(t)}$ is -- at a glance -- not necessarily well-defined for all $t\in[0,T]$. However, it is shown in \cite[Lem.~4.7.9]{Staffans2005} that the second integral of $\spvek {x}{u}$ is continuous as a~mapping from $[0,T]$ to $\dom(A\&B)=\dom S$. As a~consequence, the output can -- in the distributional sense -- be defined as the second derivative of $C\&D$ applied to the second integral of $\spvek {x}{u}$. This can be used to show that $(x,u,y)$ is a~generalized trajectory for \eqref{eq:ODEnode} if, and only if, $(x,u)$ is a~generalized trajectory for \eqref{e:ODE}, and (cf.\ (4.7.6) in \cite{Staffans2005})
\[
y=\left(t\mapsto\tfrac{{\rm d}^2}{{\rm d}t^2}\,C\&D\int_0^t(t-\tau)\spvek {x(\tau)}{u(\tau)}{\rm d}\tau\right)\in L^2([0,T];\calY).
\]
\end{enumerate}
\end{rem}

Next we recall the important concept of well-posed systems.

\begin{defn}[Well-posed systems]\label{def:wp}
Let $S = \sbvek{A\& B}{C\& D}$ be a~system node  on
$(\calY,\calX,\calU)$. The system \eqref{eq:ODEnode} is called
\emph{well-posed}, if for some (and hence all) $T>0$, there exists
some $c_T>0$, such that the classical (and thus also the generalized)
trajectories for \eqref{eq:ODEnode} on $[0,T]$ fulfill
\[
\|x(T)\|_\calX+\|y\|_{L^2([0,T];\calY)}\,\leq\,
c_T\big(\|x(0)\|_\calX+\|u\|_{L^2([0,T];\calU)}\big).
\]
\end{defn}

\begin{rem}[Well-posed systems]\label{rem:wp}
Let $S = \sbvek{A\& B}{C\& D}$ be a~system node on
$(\calY,\calX,\calU)$ and $T>0$. Well-posedness of \eqref{eq:ODEnode}
is equivalent to the well-definition and the boundedness of the mappings
\[
\begin{aligned}
\frakB_T\colon &&L^2([0,T];\calU)\to&\, \calX,\quad&\frakC_T\colon && \calX\to
&\, L^2([0,T];\calY),\\
\frakD_T\colon &&L^2([0,T];\calU)\to&\, L^2([0,T];\calY),
\end{aligned}
\]
defined by
\begin{itemize}
\item[--] $\frakB_T u=x(T)$, where $(x,u)$ is the generalized trajectory for \eqref{e:ODE} on $[0,T]$ with initial value $x(0)=0$;
\item[--] $\frakC_T x_0=y$, where $(y,x,u)$ is the generalized trajectory for \eqref{eq:ODEnode} on $[0,T]$ with input $u=0$ and initial value $x(0)=x_0$;
\item[--] $\frakD_T u=y$, where $(y,x,u)$ is the generalized trajectory for \eqref{eq:ODEnode} on $[0,T]$ with initial value $x(0)=0$.
\end{itemize}
In view of Remark~\ref{rem:sols}\,\eqref{rem:sols2}, we have
\[
\frakB_T u = \int_0^T \frakA_{-1}(T-\tau)Bu(\tau){\rm d}\tau\quad \forall u\in L^2([0,T];\calU).
\]
In particular, well-posedness implies that the above integral is an element of $\calX$.
A~combination of Proposition~\ref{prop:solex} with Remark~\ref{rem:sols}\,\eqref{rem:sols2} yields that $\frakA(t)x_0\in\dom A$ for all $t>0$ and $x_0\in\dom A$.
Thus, with $C$ as in Remark~\ref{rem:nodesnew}, setting $y=C\frakA(\cdot)x_0$ and $x=\frakA(\cdot)x_0$ yields a classical trajectory $(x,0,y)$ for \eqref{eq:ODEnode} on $[0,T]$ with $x(0)=x_0$. Well-posedness implies that the mapping $x_0\mapsto C\frakA(\cdot)x_0$ has an extension to a~bounded linear operator $\frakC_T\colon \calX\to L^2([0,T];\calY)$, see \cite[Thm.~4.7.14]{Staffans2005}.
\end{rem}

%
%We briefly recall a suitable solution concept for systems naturally governed by system nodes $S=\sbvek{A\&B}{C\&D}$, that is,
%\begin{align}\label{e:system}
%\begin{bmatrix}
%\dot x(t)\\ y(t)
%\end{bmatrix}
%= \begin{bmatrix}
%A\&B\\C\&D
%\end{bmatrix}\begin{bmatrix}
%x(t)\\u(t)
%\end{bmatrix}.
%\end{align}
%A {\em classical \braces{system} trajectory} for \eqref{e:system} is a triple
%$$
%(x,u,y)\,\in\,C^1(\mathbb{R}_0^+;\calX)\times C(\mathbb{R}_0^+;\calU)\times C(\mathbb{R}_0^+;\calY)
%$$
%which for all $t\geq 0$ satisfies
%$$
%\begin{bmatrix}
%x(t)\\u(t)
%\end{bmatrix}\in \dom S
%\qquad\text{and}\qquad
%\begin{bmatrix}
%\dot x(t)\\ y(t)
%\end{bmatrix}
%= \begin{bmatrix}
%A\&B\\C\&D
%\end{bmatrix}\begin{bmatrix}
%x(t)\\u(t)
%\end{bmatrix}.
%$$
%The following result ensures the existence of unique classical trajectories with suitable control functions and initial values.
%
%\begin{lem}[{\cite[Prop.~2.4]{OpmeerStaffans2019}}]
%Let $S$ be a system node on $(\calX,\calU,\calY)$. Then, for all initial values $x_0\in \calX$ and controls $u\in W^{1,2}_\text{loc}(\mathbb{R}_0^+;\calU)$ with $\sbvek{x_0}{u(0)}\in \dom S$ there is a unique classical system trajectory $(x,u,y)$ for \eqref{e:system} with $x(0)=x_0$.
%\end{lem}
%
%\FP{Auch was ueber mild solutions sagen.}
%
\subsection{Dissipative operators}

\begin{defn}[(Maximal) dissipative operator]
Let $\calX$ be a~Hilbert space. A subspace $\calZ\subset \calX\times\calX^*$ is called {\em dissipative}, if
\[
\forall\,(v,w)\in \calZ:\qquad   \Re\<v,w\>_{\calX,\calX^*}\le 0.
\]
Further, $\calZ\subset\calX\times\calX^*$ is called {\em maximal dissipative} if it is dissipative and not a proper subset of a dissipative subspace of $\calX\times\calX^*$.\\
An operator $A:\calX\supset\dom A\to\calX^*$ is called {\em (maximal) dissipative} if the graph of $A$, i.e., $\setdef{(x,Ax)}{\,x\in \dom A}\subset\calX\times\calX^*$, is (maximal) dissipative.\\
An operator $A:\calX\supset\dom A\to\calX$ is called {\em (maximal) dissipative} if $\calR_\calX A$ is (maximal) dissipative.
\end{defn}

Hence, an operator $A:\calX\supset\dom A\to\calX$ is dissipative if, and only if, $\Re\<Ax,x\>\le 0$ for all $x\in\dom A$, which matches the usual definition.\\
Let $A : \calX\supset\dom A\rightarrow \calX^*$ be a~dissipative operator. It follows from the definition of dissipativity that $\lambda\calR_\calX-A$ is injective for all $\lambda\in\C_+$. Maximal dissipativity further implies the bijectivity of $\lambda\calR_\calX-A$.

\begin{prop}[Maximal dissipative operators]\label{prop:maxdiss}
Let $\calX$ be a~Hilbert space and consider a closed and dissipative linear operator $A:\calX\supset\dom A\to\calX^*$. Then the following statements are equivalent:
\begin{enumerate}[(i)]
\item\label{prop:maxdiss1} $A$ is maximal dissipative,
\item\label{prop:maxdiss2} $\lambda \calR_\calX- A$ is surjective for some (and hence for all) $\lambda\in\C_+$,
\item\label{prop:maxdiss3} $\lambda \calR_\calX- A$ has dense range for some (and hence for all) $\lambda\in\C_+$,
\item\label{prop:maxdiss4} $A$ is densely defined, and $\lambda \calR_\calX- A^*$ is injective for some (and hence for all) $\lambda\in\C_+$.
\item\label{prop:maxdiss5} $A$ is densely defined, and $A^*$ is dissipative (and thus also maximal dissipative).
\end{enumerate}
\end{prop}
\begin{proof}
By using that $\calR_\calX$ is the Riesz isomorphism, it suffices to consider the case where $\calX^*$ is identified with $\calX$ (consequently, $\calR_\calX=\Id_\calX$).\\
The equivalence between \eqref{prop:maxdiss1}, \eqref{prop:maxdiss2} follows from \cite[Thm.~1.6.4]{BehrHassdeSn20}, which is shown in the~larger context of linear relations. Moreover, for $x\in\dom A$ and $\la\in\C^+$, we have
\[
\|(\lambda \calR_\calX- A)x\|_{\calX^*}\|x\|_\calX\,\ge\,\Re\langle x,(\lambda \calR_\calX- A)x\rangle_{\calX,\calX^*}\,\ge\,(\Re\lambda)\;\|x\|_{\calX}^2.
\]
This shows that $\lambda\calR_\calX - A$ has per se closed range for all $\la\in\C^+$. Hence, \eqref{prop:maxdiss2} and \eqref{prop:maxdiss3} are equivalent.\\
Moreover, an application of \cite[Prop.~3.1.6]{TuWe09} yields that maximal dissipative operators are densely defined. Then, an application of
\cite[Thm.~1.6.4 \& Prop.~1.6.7]{BehrHassdeSn20} further shows that each of \eqref{prop:maxdiss4} and \eqref{prop:maxdiss5} is equivalent to $A$ being maximal dissipative.
%Our statement now follows from the fact that $A:\calX\supset\dom A\to\calX^*$ is maximal dissipative if, and only if,
% $\calR_\calX^{-1}A:\calX\supset\dom A\to\calX$ is maximal dissipative.
\end{proof}

\subsection{Positive sesquilinear forms and quasi Gelfand triples}\label{sec:gelf}

Sesqulinear forms are -- as inner products -- assumed to be linear in the first argument and antilinear in the second one. In the following, we mainly consider sesquilinear forms with additional properties.

\begin{defn}[Positive/closed sesquilinear forms]
Let $\calX$ be a~Hilbert space. A~sesquilinear form $h:\dom h\times\dom h\to\C$ for some dense subspace $\dom h\subset\calX$ is called positive, if
\[
\forall\, x\in\dom h\setminus\{0\}:\quad h(x,x)>0.
\]
A positive sesquilinear form $h:\dom h\times\dom h\to\C$ is called {\em closed}, if $\dom h$ is complete with respect to the norm
\begin{equation}\|x\|_{\dom h}:=\big(\|x\|_{\calX}^2+h(x,x)\big)^{1/2}.\label{eq:normdomh}
\end{equation}
\end{defn}

For a~sesqulinear form $h$, realness of $h(x,x)$ for all $x\in\dom h$ implies that $h$ is symmetric in the sense that $h(x,y)=\overline{h(y,x)}$ for all $x,y\in\dom h$ \cite[Sec.~6.1]{Kato2013}. As a~consequence,
%$x\mapsto \big(h(x,x)\big)^{1/2}$ is a~seminorm, if $h$ is nonnegative. If $h$ is positive, then it
$h$ defines an inner product on $\dom h$, but this does not necessarily yield a~Hilbert space.

%For a~nonnegative sesquilinear form $h:\dom h\times\dom h\to\C$, we set
%\[\ker h=\setdef{x\in\dom h}{h(x,x)=0}.\]
%It can be easily verified that $\ker h$ is a~subspace of $\dom h$. Closedness of $h$ further implies that $\ker h$ is a~closed subspace of $\calX$.
%\begin{prop}\label{prop:h_P}
%Let $\calX$ be a~Hilbert space, and let $h:\dom h\times\dom h\to\C$ be a~densely defined, nonnegative and closed sesquilinear form. Further, let $P\in L(\calX)$ be the orthogonal projector onto $\ker h$ (with respect to the inner product in $\calX$). Then the %sesquilinear form
%\begin{equation}\label{eq:h_P}
%\begin{aligned}
%h_P:\qquad\dom h\times\dom h&\to \C,\\
%(x,y)&\mapsto h(x,y)+\langle x,Py\rangle_{\calX}
%\end{aligned}
%\end{equation}
%is positive and closed.
%\end{prop}
%\begin{proof}
%Positivity if $h_P$ is obvious. It remains to prove closedness. Define, for $x\in\dom h$,
%\begin{align*}
%\|x\|_{h}&:=\big(\|x\|_{\calX}^2+h(x,x)\big)^{1/2},\\
%\|x\|_{h,P}&:=\big(\|x\|_{\calX}^2+h_P(x,x)\big)^{1/2}.
%\end{align*}
%Then, for all $x\in \dom h$, we have
%\begin{multline*}\|x\|_{h}^2=\|x\|_{\calX}^2+h(x,x)\leq\underbrace{\|x\|_{\calX}^2+\|Px\|_{\calX}^2+h(x,x)}_{=\|x\|_{h,P}^2}\leq 2
%\big(\|x\|_{\calX}^2+h(x,x)\big)=2\,\|x\|_{h}^2.\end{multline*}
%Hence, the norms $\|\cdot\|_{h}$ and $\|\cdot\|_{h,P}$ are equivalent. Since closedness of $h$ yields that $\dom h$ is complete with respect to
%$\|\cdot\|_{h}$, the above proven norm equivalence shows that $\dom h$ is complete with respect to
%$\|\cdot\|_{h,P}$. This shows that $h_P$ is a~closed form.
%\end{proof}
The concept of quasi Gelfand triple from \cite{Skrepek21} plays a~fundamental role for our considerations on Hamiltonians, see also the recent preprint \cite{Skr23}.

\begin{defn}[Quasi Gelfand triple]\label{def:qgt}
Let $\calX$ be a~Hilbert space with inner product $\langle\cdot,\cdot\rangle_{\calX}$, and let $h:\dom h\times \dom h\to\C$ be a~densely defined, closed and positive sesquilinear form. Let $\calX_h$ be the completion of $\dom h$
with respect to the norm $\|x\|_h=h(x,x)^{1/2}$. Further, consider %the anti-dual $\calX_{-1}$ of $\calX_1$ with respect to the pivot space $\calX$. That is,
\[\|x'\|_{h-}=\sup_{x\in\dom h\setminus\{0\}}\frac{|\langle x',x\rangle_{\calX}|}{\|x\|_{h}},\quad \tilde{\calX}_{h-}=\setdef{x\in\calX}{\|x\|_{h-}<\infty},\]
and let $\calX_{h-}$ be the completion of $\tilde{\calX}_{h-}$ with respect to the norm $\|\cdot\|_{h-}$.\\
Then we call $(\calX_{h-},\calX,\calX_{h})$ the {\em quasi Gelfand triple associated with $h$}.
\end{defn}

\begin{rem}[Quasi Gelfand triples]\label{rem:gelftrip}
Let $\calX$ be a~Hilbert space, let $h:\dom h\times \dom h\to\C$ be a~closed positive sesquilinear form, and let $(\calX_{h-},\calX,\calX_{h})$ be the quasi Gelfand triple associated with $h$.
\begin{enumerate}[(a)]
\item\label{rem:gelftrip0} It is shown in \cite[Cor.~4.8]{Skrepek21} that $\calX_{h-}$ is, in a~canonical way, isometrically isomorphic to $\calX_{h}^*$.
    It is also referred to as {\em the anti-dual of $\calX_{h}$ with respect to the pivot space $\calX$}.\\
    We canonically identify $\calX_{h}^*=\calX_{h-}$, and therefore we will speak of the quasi Gelfand triple $(\calX_{h}^*,\calX,\calX_{h})$.
\item\label{rem:gelftrip1} If $h$ is bounded in the sense that there exists some $c>0$ with $h(x,x)\leq c\,\|x\|_\calX^2$, then closedness of $h$ gives $\calX\subset\calX_h$. The construction of the quasi Gelfand triple then yields $\calX_h^*\subset\calX\subset\calX_h$.
\item\label{rem:gelftrip2} If $h$ is coercive in the sense that there exists some $c>0$ with $h(x,x)\geq c\,\|x\|_\calX^2$, then we obtain $\calX_h=\dom h$, whence, in particular, $\calX_h\subset\calX$. The construction of the quasi Gelfand triple then yields $\calX_h\subset\calX\subset\calX_h^*$.
\item A quasi Gelfand triple with $\calX_h^*\subset\calX\subset\calX_h$ or $\calX_h\subset\calX\subset\calX_h^*$ is called {\em Gelfand triple}.
\item \label{rem:gelftrip3} $\calX\cap\calX_h\cap\calX_{h}^*$ is dense in all the spaces $\calX$, $\calX_h$ and $\calX_{h}^*$ with their respective norms \cite[Prop.~4.14 \& Cor.~4.16]{Skrepek21}.
\end{enumerate}
\end{rem}

Assume that $h:\dom h\times\dom h\to\C$ is a~densely defined and closed positive sesquilinear form on $\calX$. Due to Kato's second representation theorem \cite[Chap.~6.2]{Kato2013}, there exists a unique positive self-adjoint operator $H:\calX\supset\dom H\to\calX$, such that its operator square root fulfills $\dom H^{1/2} = \dom h$, and
\begin{align}\label{eq:Hop}
\forall\,x,y\in\dom h:\quad
h(x,y) = \big\langle H^{1/2}x,H^{1/2}y\big\rangle.
\end{align}
The domain of the operator $H$ is moreover dense in $\dom h$, equipped with the norm \eqref{eq:normdomh}. We call $H$ the {\em operator associated with $h$}.
%Note that positivity of $h$ is equivalent to $\ker H = \{0\}$.
%Clearly, we can do without the scaling factor $\tfrac12$. This factor is however notationally advantageous in conjunction with Hamiltonians, as we will see in the forthcoming section.

We further collect some connections between densely defined and closed positive sesquilinear forms, the operator $H$ as in \eqref{eq:Hop},
and quasi Gelfand triples associated with $h$.

\begin{prop}[Quasi Gelfand triples]\label{prop:sesqRiesz}
Let $\calX$ be a~Hilbert space, let $h:\dom h\times \dom h\to\C$ be a~closed positive sesquilinear form, and let $(\calX_{h}^*,\calX,\calX_{h})$ be the quasi Gelfand triple associated with $h$. Then the following statements hold:
\begin{enumerate}[(a)]
\item\label{prop:sesqRiesz0} $\ran H^{1/2}\subset\calX_h^*$.
\item\label{prop:sesqRiesz1} The restriction of the inner product $\langle\cdot,\cdot\rangle_{\calX}$ to $\big(\calX_{h}^*\cap\calX\big)\times \big(\calX_{h}\cap\calX\big)$ extends to the duality product $\langle\cdot,\cdot\rangle_{\calX_h^*,\calX_h}$.
\item\label{prop:sesqRiesz2} $H^{1/2}:\calX\supset\dom h\to \calX$ uniquely extends to a~bounded linear operator $U:\calX_h\to \calX$. This operator is moreover an isometric isomorphism.
\item\label{prop:sesqRiesz3} $H^{1/2}:\calX\supset\dom h\to \calX$ uniquely extends to a~bounded linear operator $V:\calX\to \calX_h^*$. This operator fulfills $V=U^*$. In particular, it is an isometric isomorphism.
\item\label{prop:sesqRiesz4} $H:\dom H\to \calX$ uniquely extends to a~bounded linear operator $\widetilde{H}:\calX_h\to \calX_h^*$. This operator coincides with the Riesz isomorphism, i.e., $\widetilde{H}=\calR_{\calX_h}$.
\end{enumerate}
\end{prop}
\begin{proof}\
\begin{enumerate}[(a)]
\item For $x,z\in\dom h\backslash\{0\}$ we have
\[
\frac{|\<H^{1/2}z,x\>_\calX|}{\|x\|_h} = \frac{|\<z,H^{1/2}x\>_\calX|}{\|H^{1/2}x\|_\calX}\,\le\,\|z\|_\calX,
\]
so that, indeed, $H^{1/2}z\in\tilde\calX_{h-}\subset\calX_{h-} = \calX_h^*$.
\item This is shown in \cite[Rem.~4.11]{Skrepek21}.
\item By \eqref{eq:Hop}, for $x\in\dom h$  we have $\|x\|_h^2 = \|H^{1/2}x\|_{\calX}^2$. Moreover, by construction of $\calX_h$, $\dom h$ is dense in
$\calX_h$. Hence, $H^{1/2}$ extends uniquely to an isometric linear
operator $U:\calX_h\to \calX$.\\ %Further, since $H^{1/2}:\calX\supset\dom
%h\to\calX$ is self-adjoint, we have
%\[\forall\,x,y\in\dom h:\quad
%\langle x,H^{1/2}y\rangle_{\calX}=\langle
%H^{1/2}x,y\rangle_{\calX}.\]
To complete the proof, it remains to show that $U$ is onto: Since
$U$ is norm-preserving, $\ran U$ is closed. Hence, it suffices to
prove that $(\ran U)^\bot=\{0\}$. Assume that $z\in(\ran U)^\bot$.
Then
\begin{equation}
\forall \,x\in\dom h:\quad 0=\langle z,Ux\rangle_{\calX}
=\langle z,H^{1/2}x\rangle_{\calX}.\label{eq:h0}
\end{equation}
Hence, $z\in (\ran H^{1/2})^\perp = \ker H^{1/2} = \{0\}$.
%Therefore, $z\in \dom \big(H^{1/2}\big)^*$, and, by using that
%$H^{1/2}$ is self-adjoint, we are led to $z\in \dom H^{1/2}=\dom h$. Now let $(z_n)$ be
%a~sequence in $\dom H$ which converges in the norm
%\eqref{eq:normdomh} to $z$. By using an argumentation as in
%\eqref{eq:h0}, we obtain
%\[h(z,z)\stackrel{n\to\infty}{\longleftarrow}h(z,z_n)
%=\langle z,U \, H^{1/2}z_n\rangle_{\calX}=0,\]
%whence $h(z,z)=0$. Positivity of $h$ now gives $z=0$.
\item Consider the operator $V=U^*\in L(\calX,\calX_h^*)$. Then $V$ is an isometric isomorphism by \eqref{prop:sesqRiesz2}. We further conclude from \eqref{prop:sesqRiesz0} and \eqref{prop:sesqRiesz1} that for all $x,y\in\dom h$ we have
\[
\langle Vx,y\rangle_{\calX_h^*,\calX_h}=\langle x,Uy\rangle_{\calX}=\langle x,H^{1/2}y\rangle_{\calX}=\langle H^{1/2}x,y\rangle_{\calX} = \<H^{1/2}x,y\>_{\calX_h^*,\calX_h}.
\]
Density of $\dom h$ in $\calX_h$ yields that $V$ extends $H^{1/2}$.
%, i.e.,
%\[
%\forall x\in\dom h:\quad Vx=H^{1/2}x.
%\]
Another use of the fact that $\dom h$ is dense in $\calX$ further implies that this extension is unique.
\item Consider the operator $\widetilde{H}=VU\in L(\calX_h,\calX_h^*)$, which is an isometric isomorphism by a~combination of \eqref{prop:sesqRiesz2} and \eqref{prop:sesqRiesz3}. For $x\in\dom H$, we have
\[
Hx=H^{1/2}H^{1/2}x=VUx=\widetilde{H}x,
\]
whence $\widetilde{H}$ is truly an extension of $H$. By density of $\dom H$ in $\calX_h$, this extension is further unique. It remains to prove that $\widetilde{H}=\calR_{\calX_h}$. By using that $V=U^*$ together with the fact that $U:\calX_h\to\calX$ is an isometric isomorphism, the desired result follows from
\[
\forall x,y\in\calX_h:\quad\langle \widetilde{H}x,y \rangle_{\calX_h^*,\calX_h}=
\langle VUx,y \rangle_{\calX_h^*,\calX_h}
=\langle Ux,Uy \rangle_{\calX}=\langle x,y \rangle_{\calX_h}.
\]
\end{enumerate}
\end{proof}

\section{Port-Hamiltonian system nodes
% with non-degenerate Hamiltonian
}\label{sec:phsysnodes}

In this main part of the article, we introduce our approach to infinite-dimensional port-Hamiltonian systems. It is divided into three parts: First, we define Hamiltonians by means of densely defined and closed positive sesquilinear forms. Thereafter, we introduce the concept of {\em dissipation node} which is the counterpart of the matrix $ M$ in \eqref{e:phs_compact2}. Hamiltonians and dissipation nodes are the ingredients for port-Hamiltonian systems, which are defined in the last part of this chapter.

\subsection{Quadratic Hamiltonians and the energetic space}\label{sec:quham}

As this article is devoted to linear problems, we restrict to quadratic Hamiltonians (which are referred to as ``Hamiltonians'' for the sake of brevity) throughout this work. These are associated with positive sesquilinear forms.

\begin{defn}[Hamiltonian]
Let $\calX$ be a~Hilbert space, and let $h:\dom h\times\dom h\to\C$ be
a~densely defined and closed positive sesquilinear form. Then the mapping
\begin{align*}
\mathcal{H}:\qquad\dom h&\to\R,\\
x&\mapsto \tfrac12\, h(x,x)
\end{align*}
is called~{\em Hamiltonian associated with $h$}.
\end{defn}

%Let $P\in L(\calX)$ be an orthogonal projector onto $\ker h$.
%Proposition~\ref{prop:h_P} yields that $h_P$ as in \eqref{eq:h_P}
%fulfills $\dom h_P=\dom h$, and it is positive and closed. We use this form for the construction of some spaces which play an extraordinary
%role throughout this work.
In the context of Hamiltonian systems, the Hamiltonian expresses the energy. This justifies the following definition.

\begin{defn}[Energetic space]\label{d:energetic}
Let $\calX$ be a~Hilbert space, and let $h:\dom h\times\dom h\to\C$ be
a~densely defined and closed positive sesquilinear form. Let $(\calX_h^*,\calX,\calX_h)$ be the quasi Gelfand triple associated with $h$. Then we call $\calX_h$ the {\em energetic space}.
\end{defn}

\begin{rem}[Energetic space]
In the case that the form $h$ is coercive the term {\em energetic space} has already been used in the literature (see, e.g., \cite{zeidler2012}). {\sc Johnson} \cite{johnson1987} calls the norm $\|\cdot\|_h$ the {\em energy norm}. The extension $\widetilde{H}:\calX_h\to\calX_h^*$ of the operator $H$ associated with $h$ is called the {\em energetic extension} of $H$. We adopt these nomenclatures and raise them to the more general case of quasi Gelfand triples.
\end{rem}

The definition of the energetic space yields that the Hamiltonian $\mathcal{H}$ extends to $\calX_h$, for which we use the same symbol. In fact, we have
%, where $\widetilde{P}\in L(\calX_h,\calX_h^*)$ is - by a density argument - uniquely defined by $\widetilde{P}x=Px$ for all $x\in\dom h$.
%As a~consequence, the Hamiltonian fulfills
\begin{equation}\label{eq:Hamnorm}
\forall\,x\in\calX_h:\quad\mathcal{H}(x) = \tfrac 12\,\<\wt Hx,x\>_{\calX_h^*,\calX_h} = \tfrac12\,\|x\|_{\calX_h}^2,
\end{equation}
where $\widetilde{H}:\calX_h\to\calX_h^*$ denotes the energetic extension of the operator $H$ associated with $h$.
%where $Q\in L(\calX_h)$ is the orthogonal projector with $\ker Q=\ker h$.
Further note that for all $x,y\in\calX_h$,
\begin{align*}
 \lim_{h\to 0}\frac{\mathcal{H}(x+hy)-\mathcal{H}(x)}{h}&=
\lim_{h\to 0}\frac{\langle \widetilde{H}(x+hy),x+hy\rangle_{\calX_h^*,\calX_h}-\langle \widetilde{H}x,x\rangle_{\calX_h^*,\calX_h}}{2h}\\
&=\lim_{h\to 0}\frac{2h\Re\langle \widetilde{H}x,y\rangle_{\calX_h^*,\calX_h}+h^2\langle \widetilde{H}y,y\rangle_{\calX_h^*,\calX_h}}{2h}
=\Re\langle \widetilde{H}x,y\rangle_{\calX_h^*,\calX_h},
\end{align*}
whence $\widetilde{H}$ is the G\^{a}teaux derivative of $\mathcal{H}:\calX_h\to\R$.

\subsection{Dissipation nodes}

This part is devoted to a~suitable generalization of $ M$ in \eqref{e:phs_compact} to the infinite-dimensional situation. As in the part on system nodes, we denote the canonical projection onto $\calX_h$ ($\calU^*$) in $\calX_h\times\calU^*$ by $P_{\calX_h}$ ($P_{\calU^*}$, respectively). Further, we note that, by the fact that Hilbert spaces are reflexive, we can canonically identify $(\calX_h^*\times\calU)^*=\calX_h\times\calU^*$. %Also recall that $\calR_{\calX_h}:\calX_h\to\calX_h^*$ is the Riesz map.

\begin{defn}[Dissipation node]\label{def:dissnode}
Let $\calX_h$, $\calU$ be Hilbert spaces with anti-duals $\calX_h^*$, $\calU^*$.
A {\em dissipation node} on $(\calX_h,\calU)$ is a (possibly unbounded) linear operator $ M:\calX_h^*\times\calU\supset\dom M\to\calX_h\times\calU^*$ that satisfies
\begin{enumerate}[(a)]
\item\label{def:dissnode1} $M$ is closed and dissipative;
\item\label{def:dissnode2} $P_{\calX_h} M:\calX_h^*\times\calU\supset\dom M\to\calX_h$ is closed;
\item\label{def:dissnode3} for all $u\in \calU$, there exists some $x'\in \calX_h^*$ with $\spvek{x'}{u}\in \dom M$;
\item\label{def:dissnode4} for the {\em main operator} $F : \calX_h^*\supset\dom F\to\calX_h$ with
\[
\dom F := \setdef{x'\in\calX_h^*}{\spvek {x'}0\in\dom  M}
\]
and $Fx' := P_{\calX_h}  M\spvek {x'}0$, there exists some $\lambda\in\C_+$ such that $\lambda \calR_{\calX_h}^{-1}-F$ has dense range.
\end{enumerate}
\end{defn}

We set
\[F\& G := P_{\calX_h}  M,
\qquad\qquad
K\& L := P_{\calU^*}  M,
\]
so that $M$ can be written as
\[
 M = \sbvek{F\& G}{K\& L}.
\]
We collect some properties of dissipation nodes.

\begin{prop}[Dissipation nodes]\label{prop:dissnode0}
Let $\calX_h$, $\calU$ be Hilbert spaces with anti-duals $\calX_h^*$, $\calU^*$, and let
$M=\sbvek{F\&G}{K\&L}:\calX_h^*\times\calU\supset\dom M\to\calX_h\times\calU^*$ be a~dissipation node on $(\calX_h,\calU)$. Then the following statements hold:
\begin{enumerate}[(a)]
\item\label{prop:dissnode02} $F$ is maximal dissipative.
\item $F$ and $M$ are densely defined.
\item\label{prop:dissnode00} For all $\lambda\in\C_+$, $\lambda\calR_{\calX_h}^{-1}-F:\calX^*_h\supset\dom F\to\calX_h$ is bijective. In particular, $\lambda\calR_{\calX_h}^{-1}-F$ has a~bounded inverse for all $\lambda\in\C_+$.
\item\label{prop:dissnode04} $F\calR_{\calX_h}$ generates a~contractive semigroup $\frakA(\cdot)$ on $\calX_h$. That is, $\frakA(\cdot)$ is a~strongly continuous semigroup on $\calX_h$ with $\|\frakA(t)\|_{L(\calX_h)}\leq1$ for all $t\ge0$.
\end{enumerate}
\end{prop}
\begin{proof}\
\begin{enumerate}[(a)]
\item First, note that $F$ is closed by property \eqref{def:dissnode2} in Definition~\ref{def:dissnode}. The dissipativity of $F$ follows directly from the dissipativity of $M$. Moreover, property \eqref{def:dissnode4} in Definition~\ref{def:dissnode} together with Proposition~\ref{prop:maxdiss} yields that $M$ is maximal dissipative.

\item Combining \eqref{prop:dissnode02} with Proposition~\ref{prop:maxdiss}, we obtain that $\dom F$ is dense in $\calX_h$.
Density of $\dom M$ in $\calX^*_h\times\calU$ is now a~consequence of density of $\dom F$ in $\calX_h^*$ together with property \eqref{def:dissnode3} in Definition~\ref{def:dissnode}.

\item The first statement follows from Proposition~\ref{prop:maxdiss}, the second one from closedness of $F$.

\item By using that $F$ is maximal dissipative, we have that $F\calR_{\calX_h}:\calX_h\supset\dom(F\calR_{\calX_h})\to \calX_h$ with
$\dom (F\calR_{\calX_h})=\calR_{\calX_h}^{-1}\dom F$ is again maximal dissipative. Then the Lumer-Philips theorem
\cite[Prop.~3.8.4]{TuWe09} implies that $F\calR_{\calX_h}$ generates a~contractive semigroup on $\calX_h$.
\end{enumerate}
\end{proof}

Our definition of dissipation node strongly resembles that of operator nodes. The connection between those concepts is highlighted in the following result.

\begin{prop}[Dissipation nodes and dissipative system nodes]\label{prop:dissnode}
Let $\calX_h$, $\calU_h$ be Hilbert spaces with anti-duals $\calX_h^*$, $\calU_h^*$, and let $ M:\calX_h^*\times\calU\subset\dom M\to \calX_h\times\calU^*$ be a~linear operator. Then the following statements are equivalent:
\begin{enumerate}[(i)]
\item\label{prop:dissnode1} $M$ is a ~dissipation node on $(\calX_h,\calU)$.
\item\label{prop:dissnode2} $N:= M\sbmat{\calR_{\calX_h}}{0}{0}{\Id_{\calU}}:\calX_h\times \calU\supset\dom M\to\calX_h\times \calU^*$ with
$\dom N = \sbmat{\calR_{\calX_h}^{-1}}{0}{0}{\Id_{\calU}}\dom M$
 is a~dissipative system node on $(\calU^*,\calX,\calU)$.
\end{enumerate}
\end{prop}
\begin{proof}\

\noindent``\eqref{prop:dissnode1}$\Rightarrow$\eqref{prop:dissnode2}'':
Assume that $M$ is a ~dissipation node on $(\calX_h,\calU)$. Then $M$ is dissipative, and it has properties \eqref{def:dissnode1}, \eqref{def:dissnode2} and \eqref{def:dissnode3} in Definition~\ref{def:dissnode}. This yields that $N$ fulfills \eqref{def:opnode1}, \eqref{def:opnode2} and \eqref{def:opnode3} in Definition~\ref{def:opnode}. Further, by Proposition~\ref{prop:dissnode0}, $F\calR_{\calX_h}$ generates a~strongly continuous semigroup on $\calX_h$. Altogether, we obtain that $N$ is a~dissipative system node.\\
``\eqref{prop:dissnode2}$\Rightarrow$\eqref{prop:dissnode1}'':
Assume that $N$ is a~dissipative system node on $(\calU^*,\calX,\calU)$. Again, by comparing Definitions~\ref{def:opnode}\&\ref{def:dissnode}, we directly obtain that $M$ has properties \eqref{def:dissnode1}, \eqref{def:dissnode2} and \eqref{def:dissnode3} in Definition~\ref{def:dissnode}. Dissipativity of $N$ implies that $M$ has this property, too. Consequently, $F$ is dissipative. Moreover, since $F\calR_{\calX_h}$ generates a~strongly continuous semigroup on $\calX_h$, the Hille-Yosida theorem \cite[Sec.~II.3]{engel2000one} yields that $\lambda \Id_\calX-F\calR_{\calX_h}$ is boundedly invertible for some $\lambda\in\C_+$. This shows that $M$ is a~dissipation node on $(\calX_h,\calU)$.
\end{proof}

\begin{rem}[Dissipation nodes]\label{rem:dissnode} \
\begin{enumerate}[(a)]
\item\label{rem:dissnode0} Analogous to Remark~\ref{rem:nodesnew}, a~dissipation node $ M = \sbvek{F\&G}{K\& L}$ on $(\calX_h,\calU)$ fulfills, by using
Lemma~\ref{lem:closedblock}, $K\&L\in \mathcal L(\dom(F\&G),\calU^*)$, where $\dom(F\&G)$ is equipped with the graph norm of $F\&G$.\\
On the other hand, if $ M = \sbvek{F\&G}{K\& L}$ is dissipative and has properties \eqref{def:dissnode2}--\eqref{def:dissnode4} in Definition~\ref{def:dissnode} with, additionally, $K\&L\in \mathcal L(\dom(F\&G),\calU^*)$, then Lemma~\ref{lem:closedblock} yields that $ M$ is a~dissipation node.
\item\label{rem:dissnode1} As for system nodes, we can define $K:\dom F\to\calU^*$ by $Kx'=K\&L\spvek{x'}0$. The previous remark yields $K\in \mathcal L(\dom F,\calU^*)$ in case where $ M = \sbvek{F\&G}{K\& L}$ is a~dissipation node.
\item\label{rem:dissnode2} {Take $\alpha\in \mathbb{C}$ such that $\alpha \calR_{\calX_h}^{-1}-F$ is bijective} and let $\calX_{h,-1}$ be the completion of $\calX_h$ with respect to the norm \[\|x\|_{\calX_{h,-1}}:= \|(\alpha \calR_{\calX_h}^{-1}-F)^{-1}x\|.\]
    Since $\calR_{\calX_h}:\calX_h\to\calX_h^*$ is an isometry, and $\alpha \calR_{\calX_h}^{-1}-F=(\alpha \Id-F\calR_{\calX_h})\calR_{\calX_h}^{-1}$, we can conclude from \cite[Prop.~2.10.2]{TuWe09} that the topology of $\calX_{h,-1}$ does not depend on the particular choice of $\alpha\in\C$. Further, $F$ and $F\& G$ extend to bounded operators $F_{-1}:\calX_{h}^*\to \calX_{h,-1}$ and $[F_{-1}\ G]\colon \calX_h^*\times \calU\to \calX_{h,-1}$, respectively.
%\item\label{rem:dissnode3} By Proposition~\ref{prop:maxdiss}, maximal dissipativity of $F$ implies property (d) in Definition~\ref{def:dissnode}. Likewise, dissipativity of $ M$ together with dissipativity of $F^*$ also yields that (d) holds.
\item\label{rem:dissnode4} Assume that $h$ is bounded and coercive. Then the inner products in $\calX_h$ and $\calX_h^*$ are both equivalent to that in $\calX$, which yields $\calX=\calX_h=\calX_h^*$. Moreover, the construction of quasi Gelfand triples yields that the duality product equals the inner product in $\calX$. Hence, in this case, $ M$ is a~dissipation node on $(\calX_h,\calU)$ if, and only if, it is a~dissipation node on $(\calX,\calU)$.
\end{enumerate}
\end{rem}

Now we are able to show that dissipation nodes have some further properties.

\begin{prop}[Dissipation nodes - further properties]\label{prop:dissnode_2}
Let $\calX_h$, $\calU$ be Hilbert spaces with anti-duals $\calX_h^*$, $\calU^*$, and let $M = \sbvek{F\&G}{K\&L}:\calX_h^*\times\calU\supset\dom  M\to\calX_h\times\calU^*$ be a~dissipation node on $(\calX_h,\calU)$. Then the following statements hold:
\begin{enumerate}[(a)]
\item\label{prop:dissnode20} $ M$ is maximal dissipative.
\item\label{prop:dissnode21} $ M^*:\calX_h^*\times\calU\supset\dom  M^*\to\calX_h\times\calU^*$ is a~dissipation node on $(\calX_h,\calU)$ with
\[
M^*=\left[\begin{smallmatrix}{[F\& G]^d}\\{[K\& L]^d}\end{smallmatrix}\right]
\]
and domain
\[\dom  M^*=\setdef{\spvek {x'}y\in \calX_h\times\calU}{F^*_{-1}x+K^*x'\in \calX},\]
where $[F\&G]^d$ is the restriction of $[F^*_{-1}\; K^*]$ to $\dom  M^*$, and, for all $\alpha\in\rho(F\calR_{\calX_h})$,
\[
[K\&L]^d\spvek {x'}y:=G^*\Big(x'-(\overline{\alpha}\calR_{\calX_h}^{-1}-F^*_{-1}\big)^{-1}K^*\Big)
+\calF(\alpha)^*y,
\]
where $\calF$ is transfer function associated with the operator node $N= M\sbmat{\calR_{\calX_h}}{0}{0}{\Id_{\calU}}$.\\
In particular, the main operator of $ M^*$ is $F^*$.
%\item\label{prop:dissnode25} $N:= M\sbmat{\calR_{\calX_h}}{0}{0}{\Id_{\calU}}:\calX_h\times \calU\supset\dom M\to\calX_h\times \calU^*$ is a~system node on $(\calU^*,\calX_h,\calU)$.
\end{enumerate}
\end{prop}
\begin{proof}\
\begin{enumerate}[(a)]
\item By using Proposition~\ref{prop:dissnode}, we can infer from \cite[Lem.~4.3]{staffans2002passive} that $N:= M\sbmat{\calR_{\calX_h}}{0}{0}{\Id_{\calU}}$ is maximal dissipative. Hence, $M$ is as well maximal dissipative.
\item Since, by the previous statement, $M$ is maximal dissipative, $M^*$ is dissipative by Proposition~\ref{prop:maxdiss}. The remaining statements follow by a combination of Proposition~\ref{prop:dissnode} with Remark~\ref{rem:adjnode}.
%\item Maximal dissipativity of $ M$ gives rise to dissipativity of $ M^*$. Therefore, its main operator (which is $F^*$ by \eqref{prop:dissnode21}) is dissipative, which implies injectivity of $\lambda \calR_{\calX_h}^{-1}-F^*$ for all $\lambda\in\C_+$. Then an application of Proposition~\ref{prop:maxdiss} yields that $F$ is maximal dissipative.
%\item This is a~direct consequence of \eqref{prop:dissnode22} and Proposition~\ref{prop:maxdiss}.
%\item This can be concluded by combining \eqref{prop:dissnode24} with Proposition~\ref{prop:dissnode}.
\end{enumerate}
\end{proof}

We further present a~result on a~certain redefinition of the second block component of a~dissipation {node}. In applications, this for instance refers to setting certain input and output channels as inactive.

\begin{prop}[Dissipation nodes and input/output redefinition]\label{prop:iored}
Let $\calX_h$, $\calU$, $\calU_{\rm new}$ be Hilbert spaces with anti-duals $\calX_h^*$, $\calU^*$, $\calU^*_{\rm new}$, let
$ M=\sbvek{F\&G}{K\&L}:\calX_h^*\times\calU\supset\dom  M\to\calX_h\times\calU^*$ be a~dissipation node on $(\calX_h,\calU)$, and let $U\in L(\calU_{\rm new},\calU)$. Then
\begin{align*}
&M_U:=\sbmat{\Id_{\calX_h}}00{U^*}M \sbmat{\Id_{\calX_h^*}}00{U}=\sbvek{F\&(GU)}{U^*(K\&(LU))},\\
&\dom M_U=\setdef{\spvek{x'}u}{\spvek{x'}{Uu}\in\dom M}
\end{align*}
is a~{dissipation} node on $(\calX_h,\calU^*_{\rm new})$.
\end{prop}
\begin{proof}
We successively show the properties \eqref{def:dissnode1}--\eqref{def:dissnode4} in Definition~\ref{def:dissnode} in reverse order.
\begin{itemize}
\item[\eqref{def:dissnode4}] This property is trivially fulfilled, since $M$ and $M_U$ have the same main operators.
\item[\eqref{def:dissnode3}] Let $u_{\rm new}\in \calU_{\rm new}$. Then there exists some $x'\in \calX_h^*$ with $\spvek{x'}{Uu_{\rm new}}\in \dom M$, and the definition of $M_U$ leads to $\spvek{x'}{u_{\rm new}}\in \dom M_U$.
\item[\eqref{def:dissnode2}] Consider a~sequence in $\dom M_U$ with
\begin{align*}
\left(\spvek{x_n'}{u_{{\rm new},n}}\right)&\to \spvek{x'}{u_{\rm new}}\text{ in $\calX_h^*\times\calU_{\rm new}$},\\
\left(F\&(GU)\spvek{x_n'}{u_{{\rm new},n}}\right)&\to x\text{ in $\calX_h$.}
\end{align*}
Then, by boundedness of $U$,
\begin{align*}
\left(\spvek{x_n'}{Uu_{{\rm new},n}}\right)&\to \spvek{x'}{U{u_{\rm new}}}\text{ in $\calX_h^*\times\calU$},\\
\left(F\&G\spvek{x_n'}{Uu_{{\rm new},n}}\right)&\to x\text{ in $\calX_h$.}
\end{align*}
Closedness of $F\&G$ with $\dom (F\&G)=\dom M$ now leads to $\spvek{x'}{Uu_{\rm new}}\in\dom M$ with \[\big(F\&G\spvek{x'}{Uu_{\rm new}}\big)=x.\]
In other words, $\spvek{x'}{u_{\rm new}}\in\dom M_U$ with
$\big(F\&(GU)\spvek{x'}{u_{\rm new}}\big)=x$.
\item[\eqref{def:dissnode1}] Dissipativity of $M_U$ follows immediately from dissipativity of $M$. It remains to show that $M_U$ is closed.
%Consider a~sequence with
%\begin{align*}
%\left(\spvek{x_n'}{u_{{\rm new},n}}\right)&\to \spvek{x'}{u_{\rm new}}\text{ in $\calX_h^*\times\calU_{\rm new}$},\\
%\left(M_U\spvek{x_n'}{u_{{\rm new},n}}\right)&\to \spvek{x}{y}\text{ in $\calX_h\times\calU_{\rm new}^*$.}
%\end{align*}
%Then, by (ii), $\spvek{x'}{u_{\rm new}}\in\dom M_U$ with $F\&G\spvek{x'}{Uu_{\rm new}}=x$.
By Lemma~\ref{lem:closedblock} and the fact that $M$ is a~dissipation node, there exists some $C>0$ with
\[\|K\&L\spvek{x'}{Uu}\|_{\calU^*}\leq C\,\big(\|\spvek{x'}{Uu}\|_{\calX_h^*\times\calU}+\|F\&G\spvek{x'}{Uu}\|_{\calX_h}\big)\quad\forall
\spvek{x'}{u}\in\dom  M_U.\]
This means that, for all $\spvek{x'}{u}\in\dom  M_U$,
\begin{align*}
\|U^*(K\&(LU))\spvek{x'}{u}\|_{\calU^*}
&=\|U\|_{L(\calU_{\rm new},\calU)} \|K\&(LU)\spvek{x'}{u}\|_{\calU^*}\\
&\leq C\,\|U\|_{L(\calU_{\rm new},\calU)} \,\big(\|\spvek{x'}{Uu}\|_{\calX_h^*\times\calU}+\|F\&G\spvek{x'}{Uu}\|_{\calX_h}\big)\\
&= \tilde{C}\,\big(\|\spvek{x'}{u}\|_{\calX_h^*\times\calU}+\|F\&(GU)\spvek{x'}{u}\|_{\calX_h}\big)
\end{align*}
with
\[
\tilde{C}:=C\cdot \|U\|_{L(\calU_{\rm new},\calU)}\cdot\max\{\|U\|_{L(\calU_{\rm new},\calU)},1\}.
\]
Now another application of Lemma~\ref{lem:closedblock} shows that $ M_U$ is closed.
\end{itemize}
\end{proof}

\subsection{Port-Hamiltonian systems}\label{Ssec:pH}
Having introduced Hamiltonians and dissipation nodes, we are now able to introduce our concept of port-Hamiltonian systems inspired by the theory of system nodes.

\begin{defn}[Port-Hamiltonian system node]\label{def:pHnode}
Let $\calX$, $\calU$ be Hilbert spaces, let $h:\dom h\times \dom h\to\C$ with $\dom h\subset\calX$ be a~densely defined and closed positive sesquilinear form. Let $(\calX_{h}^*,\calX,\calX_{h})$ be the quasi Gelfand triple associated with $h$, and let $\widetilde{H}\in L(\calX_h,\calX_h^*)$ be the energetic extension of the operator associated with $h$ (cf.\ Proposition~\ref{prop:sesqRiesz}\,\eqref{prop:sesqRiesz4}). Further, let $ M:\calX_h^*\times\calU\supset\dom  M\to\calX_h\times\calU^*$ be a~dissipation node on $(\calX_h,\calU)$. Then we call
\begin{equation}
S:=\bmat{\Id_{\calX_h}}{}{}{-\Id_{\calU^*}} M\bmat{\widetilde{H}}{}{}{\Id_\calU}\label{eq:pHsysnode}
\end{equation}
a~{\em port-Hamiltonian system node on $(\calX_h,\calU)$}.
\end{defn}

If $M = \sbvek{F\&G}{K\&L}$, the differential equation corresponding to a~port-Hamiltonian system node is given by
\begin{equation}\label{eq:pHODEnode}
\spvek{\dot{x}(t)}{y(t)}
= \sbvek{F\&G\\[-1mm]}{-K\&L} \spvek{\widetilde{H}{x}(t)}{u(t)}.
\end{equation}
Analogous to Definition~\ref{def:traj}, a classical trajectory on $[0,T]$ is a triple
\[
(y,x,u)\,\in\, C([0,T];\calU^*)\times C^1([0,T];\calX_h)\times C([0,T];\calU)
\]
which fulfills \eqref{eq:pHODEnode} pointwise on $[0,T]$, and generalized trajectories are limits of classical trajectories in the topology of $L^2([0,T];\calU^*)\times C([0,T];\calX_h)\times L^2([0,T];\calU)$.\\
Next we collect properties of port-Hamiltonian system nodes. We show for instance that they are system nodes in the sense of Definition~\ref{def:opnode}.

\begin{thm}[Port-Hamiltonian system node]\label{thm:pHnode}
Let $S$ be a~port-Hamiltonian system node as in \eqref{eq:pHsysnode}, where we use the notation from Definition~\ref{def:pHnode}. Then the following holds:
\begin{enumerate}[(a)]
\item\label{thm:pHnode1} $S$ is a~system node on $(\calU^*,\calX_h,\calU)$;
\item\label{thm:pHnode2} the main operator $A:=F\widetilde{H}$ of ${S}$ generates a~contractive semigroup on $\calX_h$;
\item\label{thm:pHnode3} the transfer function $\calG$ associated with $S$ fulfills
\[
\forall\,s\in\C_+,\,u\in\calU:\quad\Re\langle \calG(s)u,u\rangle_{\calU^*,\calU}\ge 0;
\]
\item\label{thm:pHnode4} for all $T>0$, the generalized trajectories $(y,x,u)\,\in\,L^2([0,T];\calU^*)\times C([0,T];\calX_h)\times  L^2([0,T];\calU)$ (and thus also the classical solutions) fulfill the {\em dissipation inequality}
\begin{multline}
\forall\,t\in[0,T]:\qquad \calH(x(t))-\calH(x(0))\\=
\Re\int_0^t \left\langle  M\spvek{\widetilde{H}x(\tau)}{u(\tau)},\spvek {\widetilde{H}x(\tau)}{u(\tau)}\right\rangle_{\calX_h\times\calU^*,\calX_h^*\times\calU}{\rm d}\tau+\Re\int_0^t \langle {y(\tau)},{u(\tau)}\rangle_{\calU^*,\calU}{\rm d}\tau\\
\leq\Re\int_0^t \langle {y(\tau)},{u(\tau)}\rangle_{\calU^*,\calU}{\rm d}\tau,\label{eq:enbalinf}
\end{multline}
where $\mathcal{H}:\calX_h\to\R$ is the Hamiltonian associated with $h$.
\end{enumerate}
\end{thm}
\begin{proof}\
\begin{enumerate}[(a)]
\item Invoking that, by Proposition~\ref{prop:sesqRiesz}\,\eqref{prop:sesqRiesz4}, $\widetilde{H}=\calR_{\calX}$, we obtain from Proposition~\ref{prop:dissnode_2} that $\sbmat{\Id_\calX}{}{}{-\Id_{\calU^*}}S$ is a~system node. Then the statement follows, since the change of the sign of the operators in the output equation does not affect the property of an operator being a~system node.
\item This follows by a~combination of Proposition~\ref{prop:dissnode0}\,\eqref{prop:dissnode04} with $\widetilde{H}=\calR_{\calX}$.
\item First, we note that, by Proposition~\ref{prop:dissnode0}\,\eqref{prop:dissnode02}, $F$ is maximal dissipative, and hence so is $F\widetilde{H}=F\calR_{\calX_h}$. Thus, we conclude from Proposition~\ref{prop:maxdiss} that $s\Id-F\widetilde{H}$ is bijective for all $s\in\C_+$. Therefore, the transfer function $\calG$ is defined on $\C_+$.\\
Now let $u\in\calU$, $s\in\C_+$. By $\widetilde{H}=\calR_{\calX_h}$, we have, for $F_{-1}$ as in Remark~\ref{rem:dissnode}\,\eqref{rem:dissnode1},
\begin{align*}
\calG(s)&= -(K\widetilde{H}\&L) \,\sbvek{(s\Id-F_{-1}\widetilde{H})^{-1}G}{\Id}\\
&= -(K\calR_{\calX_h}\&L) \,\sbvek{(s\Id-F_{-1}\calR_{\calX_h})^{-1}G}{\Id}
= -K\&L \,\sbvek{(s\calR_{\calX_h}^{-1}-F_{-1})^{-1}G}{\Id}.
\end{align*}
Further,
\begin{align*}
&\phantom{=}F\&G\sbvek{(s\calR_{\calX_h}^{-1}-F_{-1})^{-1}G}{\Id}\\&=
(F-s\calR_{\calX_h}^{-1})\&G\sbvek{(s\calR_{\calX_h}^{-1}-F_{-1})^{-1}G}{\Id}+s\calR_{\calX_h}^{-1}(s\calR_{\calX_h}^{-1}-F_{-1})^{-1}G\\&=s\calR_{\calX_h}^{-1}(s\calR_{\calX_h}^{-1}-F_{-1})^{-1}G.
\end{align*}
Then, by dissipativity of $ M=\sbvek{F\&G}{K\&L}$, the latter relations imply
\begin{align*}
0&\geq \Re\left\langle\sbvek{(s\calR_{\calX_h}^{-1}-F_{-1})^{-1}G}{\Id}u,\bvek{F\&G}{K\&L}\sbvek{(s\calR_{\calX_h}^{-1}-F_{-1})^{-1}G}{\Id}u\right\rangle_{\calX_h^*\times\calU,\,\calX_h\times\calU^*}\\
&=\Re\left\langle\sbvek{(s\calR_{\calX_h}^{-1}-F_{-1})^{-1}G}{\Id}u,\spvek{s\calR_{\calX_h}^{-1}(s\calR_{\calX_h}^{-1}-F_{-1})^{-1}G}{-\calG(s)}{u}\right\rangle_{\calX_h^*\times\calU,\,\calX_h\times\calU^*} \\
&=\Re s\cdot \|(s\calR_{\calX_h}^{-1}-F_{-1})^{-1}G{u}\|_{\calX_h^*}^2-\Re\langle \calG(s)u,u\rangle_{\calU^*,\calU}\geq-\Re\langle \calG(s)u,u\rangle_{\calU^*,\calU}.
\end{align*}
\item Assume that $(y,x,u)$ is a~classical trajectory, and let $t\in[0,T]$. Then \eqref{eq:Hamnorm} together with the {product} rule gives
\begin{align*}
\tfrac{{\rm d}}{{\rm d}t}\mathcal{H}(x(t))&=
\tfrac12\tfrac{{\rm d}}{{\rm d}t}\|x(t)\|_{\calX_h}^2
=\Re\langle x(t),\dot{x}(t)\rangle_{\calX_h}^2\\
&=\Re\left\langle x(t),F\&G\spvek{\widetilde{H}x(t)}{u(t)}\right\rangle_{\calX_h}\\&=\Re\left\langle \widetilde{H}x(t),F\&G\spvek{\widetilde{H}x(t)}{u(t)}\right\rangle_{\calX_h^*,\calX_h}\\
&=\Re\left\langle \spvek{\widetilde{H}x(t)}{u(t)},\sbvek{F\&G}{K\&L}\spvek{\widetilde{H}x(t)}{u(t)}\right\rangle_{\calX_h^*\times\calU,\calX_h\times\calU^*}\\&\qquad+
\Re\left\langle u(t),-{K\&L}\spvek{\widetilde{H}x(t)}{u(t)}\right\rangle_{\calU,\calU^*}\\
&=\Re\left\langle \spvek{\widetilde{H}x(t)}{u(t)},\sbvek{F\&G}{K\&L}\spvek{\widetilde{H}x(t)}{u(t)}\right\rangle_{\calX_h^*\times\calU,\calX_h\times\calU^*}+
\Re\langle u(t),y(t)\rangle_{\calU,\calU^*}\\
&\leq \Re\langle u(t),y(t)\rangle_{\calU,\calU^*}
\end{align*}
Then \eqref{eq:enbalinf} follows by an integration of the above inequality. The case of generalized solutions can then be concluded by taking limits.
\end{enumerate}
\end{proof}

We gather some remarks on the introduced class of systems.

\begin{rem}[Port-Hamiltonian system nodes]
Let $S$ be a~port-Hamiltonian system node as in \eqref{eq:pHsysnode}, where we use the notation from Definition~\ref{def:pHnode}.
\begin{enumerate}[(a)]
\item The argumentation in the proofs of \eqref{thm:pHnode3} and \eqref{thm:pHnode4} is similar to some corresponding statements in \cite[Thm.~4.2]{staffans2002passive}
for passive systems.
\item If $\calU=\C^{m}$, then the transfer function is $\C^{m\times m}$-valued, and Theorem~\ref{thm:pHnode}~\eqref{thm:pHnode3} means that $\calG(s)+\calG(s)^*$ is positive semi-definite for all $s\in\C_+$. Transfer functions of this kind with, moreover, $\calG(s)\in\R^{m\times m}$ for all $s\in\R_{>0}$ are called {\em positive real}, and they play an outstanding role in the context of passive systems and network synthesis \cite{AndeV73}.
\item Since port-Hamiltonian system nodes are, by Theorem~\ref{thm:pHnode}~\eqref{thm:pHnode1}, system nodes in the sense of Definition~\ref{def:opnode}, we can apply Proposition~\ref{prop:solex} to see that for all $x_0\in\calX_h$ and
$u\in{W^{2,1}([0,T];\calU)}$ with $\spvek{\widetilde{H}x_0}{u(0)}\in \dom  M$ there exists a unique
classical trajectory $(y,x,u)$ for \eqref{eq:pHODEnode}.
\end{enumerate}
\end{rem}

Well-posedness in the sense of Definition~\ref{def:wp} is not a~natural property of port-Hamiltonian systems. That is, it cannot be concluded that port-Hamiltonian systems are well-posed. A~rather simple characterization for well-posedness by means of the transfer function is however possible.

\begin{prop}[Port-Hamiltonian system node and well-posedness]
Let $S$ be a~port-Hamil\-to\-nian system node as in \eqref{eq:pHsysnode}, where we use the notation from Definition~\ref{def:pHnode}. Further, let $\calG$ be the transfer function associated with $S$.
Then the system \eqref{eq:pHODEnode} is well-posed if, and only if $\calG$ is bounded on some (or equivalently, on every) vertical line in $\C_+$.
\end{prop}
\begin{proof}
By Theorem~\ref{thm:pHnode}\,\eqref{thm:pHnode1}, $S$ is a~system node on $(\calU^*,\calX_h,\calU)$. Then it follows that
\[\widetilde{S}:=\sbmat{\Id_\calX}{}{}{\calR_{\calU}^{-1}}S\] is a~system node on $(\calU,\calX_h,\calU)$, whose transfer function reads $\widetilde{\calG}(s)=\calR_{\calU}^{-1}{\calG}(s)$, and further fulfills that
$\sbmat{\Id_\calX}{}{}{-\Id}\widetilde{S}$ is dissipative. Since the Riesz map is an isometry, we have $\|\widetilde{\calG}(s)\|_{L(\calU)}=\|\calG(s)\|_{L(\calU,\calU^*)}$ and, moreover, well-posedness of  \eqref{eq:pHODEnode} is equivalent to well-posedness of the system associated with $\widetilde{S}$. By using that $S$ is a~port-Hamiltonian system node, the system associated with $\widetilde{S}$ is impedance passive in the sense of \cite{staffans2002passive}. Then the desired result can be concluded from \cite[Thm.~5.1]{staffans2002passive}.
\end{proof}

Our definition of the dissipation nodes involves the energetic space $\calX_h$, which is in turn constructed from the Hamiltonian by using the concept of quasi Gelfand triple. The ``more natural'' situation is that a~dissipation node on $(\calX,\calU)$ is given (with $\calX=\calX^*$ identified). Next we present a~criterion whether this corresponds to a~dissipation node on $(\calX_h,\calU)$ in some sense.
Such a~construction can in principle be performed as follows: We first restrict {the dissipation node $M$ on $(\calX,\calU)$} to the space
\begin{equation}
\calZ:=\setdef{\spvek{x}u\in\big(\calX_h^*\times\calU\big)\cap\dom M}{F\&G\spvek{x}u\in \calX_h},\label{eq:Zdef}\end{equation}
which is a~subspace of both $\dom M$ and {$\calX_h^*\times \calU$}. Further, the construction of $\calZ$ yields $M\calZ\subset\calX_h\times \calU^*$. That is, we obtain an operator
\begin{equation}
\widetilde{M}:\calX_h^*\times \calU\supset\dom\widetilde{M}:=\calZ\to {\calX_h\times \calU^*},\qquad \widetilde{M}\spvek{x'}u:=M\spvek{x'}u\;\forall \spvek{x'}u\in\calZ.
\label{eq:tildeMdef}\end{equation}
Thereafter, we take the closure of this operator, which is now regarded as an operator mapping from a~subspace of $\calX_h^*\times\calU$ to $\calX_h\times\calU^*$.
%To this end, we recap that, by Remark~\ref{rem:gelftrip}\,\eqref{rem:gelftrip3},
%$\calX\cap\calX_h\cap\calX_{h}^*$ is dense in all the spaces $\calX$, $\calX_h$ and $\calX_{h}^*$ with their respective norms.

%To this end, recall from Definition~\ref{def:qgt} that $\calX_h^*$ is constructed by means of the completion of a~subspace of $\calX$. This allows to form intersections of $\calX_h^*$ and subsets of $\calX$.

\begin{prop}[Dissipation nodes on $(\calX,\calU)$]\label{prop:Mclos}
Let $\calX$ and $\calU$ be Hilbert spaces and let $h:\dom h\times \dom h\to\C$ with $\dom h\subset\calX$ be a~densely defined and closed positive sesquilinear form. Let $(\calX_{h}^*,\calX,\calX_{h})$ be the quasi Gelfand triple associated with $h$. Let {$M=\spvek{F\&G}{K\&L}:\calX\times\calU\supset\dom M \to\calX\times\calU^*$} be a~dissipation node on $(\calX,\calU)$ (hereby, we identify $\calX=\calX^*$) with the following properties:
\begin{enumerate}[(a)]
%\item $\calX_h^*\times\calU\subset\dom  M$ and $\calX_h^*\times\calU\subset\dom  M^*$.
\item\label{prop:Mclosa} The space
\[
\calV:=\setdef{x\in\calX_h^*\cap\dom F}{Fx\in \calX_h},
%\calV_2&:=\setdef{x\in\calX_h^*\cap\dom F^*}{F^*x\in\dom h}
\]
is dense in $\calX_h^*$, and there exists some $\lambda\in\C_+$, such that the space
\[
\big(F-\lambda\calR_{\calX_h}^{-1}\big)\calV
\]
is dense in $\calX_h$.
\item\label{prop:Mclosb} For all $u\in\calU$, there exists some $x\in\calX_h^*$, such that $\spvek{x}{u}\in\dom  M$ with \[ M\spvek{x}{u}\in\calX_h\times\calU^*.\]
\item\label{prop:Mclosc} For all $y\in\calU$, there exists some $x\in\calX_h^*$, such that $\spvek{x}{y}\in\dom  M^*$ with
\[
M^*\spvek{x}{y}\in\calX_h\times\calU^*.
\]
\end{enumerate}
Then, for $\calZ$ as in \eqref{eq:Zdef}, the operator $\widetilde{M}$ as in \eqref{eq:tildeMdef} is closable (as an operator defined on a~subspace of $\calX_h^*\times\calU$, and mapping to $\calX_h\times\calU^*$). Its closure
\begin{equation}\label{eq:ovtildM}
\overline{\widetilde{ M}}=:\sbvek{F_c\&G_c}{K_c\&L_c}:\calX_h^*\times\calU\supset \dom\overline{\widetilde{ M}}\to\calX_h\times\calU^*.
\end{equation}
 is a~dissipation node on $(\calX_h,\calU)$.
\end{prop}
\begin{proof}~\\
{\em Step~1:} We show that $\widetilde{M}$ is closable and has a~dissipative closure: Assumptions \eqref{prop:Mclosa}\&\eqref{prop:Mclosb} yield that $\calZ$ as in \eqref{eq:Zdef} is dense in $\calX_h^*\times\calU$. Moreover, dissipativity of $M$ together with Proposition~\ref{prop:sesqRiesz}~\eqref{prop:sesqRiesz1} gives
\[
\forall\,\spvek{x}u\in\calZ:
\quad \Re\big\<\spvek{x}u,\wt M\spvek{x}u\big\>_{\calX_h^*\times\calU,\calX_h\times\calU^*} = \Re\<\spvek{x}u, M\spvek{x}u\>_{\calX\times\calU,\calX\times\calU^*}\le 0.
\]
Altogether, this means that $\widetilde{M}$ is densely defined and {dissipative}. Then \cite[Prop.~3.14~(iv)]{engel2000one} yields that it is closable with dissipative closure {$\overline{\widetilde{M}}$} .

{\em Step~2:}
We show that ${F_c\&G_c}: \calX_h^*\times\calU\supset\dom {F_c\&G_c}\to\calX_h$ with $\dom {F_c\&G_c}=\dom \overline{\widetilde{ M}}$ is closed. By
Lemma~\ref{lem:closedblock}, it suffices to show that there exists some constant $C>0$, such that
\[\|K_c\&L_c\spvek{x'}{u}\|_{\calU^*}\leq C\,\big(\|\spvek{x'}{u}\|_{\calX_h^*\times\calU}+\|F_c\&G_c\spvek{x'}{u}\|_{\calX_h}\big)\quad\forall \spvek{x'}{u}\in\dom \overline{\widetilde{ M}}.\]
Seeking for a~contradiction, we assume the converse, which means that there exists some sequence $\big(\spvek{x_n'}{u_n}\big)$ in $\calX_h^*\times\calU$, such that
\begin{align*}&(\|K_c\&L_c\spvek{x_n'}{u_n}\|_{\calU^*})\stackrel{n\to\infty}{\longrightarrow}\infty,\\ \text{ and }\;&(\|F_c\&G_c\spvek{x_n'}{u_n}\|_{\calX_h}),\; (\|\spvek{x_n'}{u_n}\|_{\calX_h^*\times \calU})\text{ bounded}.
\end{align*}
Since $\overline{\widetilde{ M}}$ is the closure of $\widetilde{ M}$, it is no loss of generality to assume that, for all $n\in\N$, $\spvek{x_n'}{u_n}\in(\calX_h^*\times\calU)\cap\dom M$ with $F\&G\spvek{x_n'}{u_n}\in\calX_h$.\\
The Banach-Steinhaus theorem \cite[Thm.~7.3]{Alt16} implies that there exists some $y\in\calU$, such that
\[\big(\langle y,K_c\&L_c\spvek{x_n'}{u_n}\rangle_{\calU,\calU^*}\big)\stackrel{n\to\infty}{\longrightarrow}\infty.\]
Further, by \eqref{prop:Mclosc}, there exists some $x\in\calX_h^*$, such that $\spvek{x}{y}\in\dom  M^*$ with $M^*\spvek{x}{y}\in\calX_h\times\calU^*$. Then Proposition~\ref{prop:sesqRiesz}~\eqref{prop:sesqRiesz1} gives
\begin{align*}
\infty\stackrel{n\to\infty}{\longleftarrow}&\phantom{=}\langle y,K_c\&L_c\spvek{x_n'}{u_n}\rangle_{\calU,\calU^*}+\langle x,F_c\&G_c\spvek{x_n'}{u_n}\rangle_{\calX_h^*,\calX_h}\\
&=\langle \spvek{x}{y},\widetilde{ M}\spvek{x_n'}{u_n}\rangle_{\calX_h^*\times\calU,\calX_h\times{\calU^*}}\\
&=\langle \spvek{x}{y},{ M}\spvek{x_n'}{u_n}\rangle_{\calX\times\calU,\calX\times\calU^*}\\
&=\langle { M}^*\spvek{x}{y},\spvek{x_n'}{u_n}\rangle_{\calX\times\calU^*,\calX\times\calU}
=\langle { M}^*\spvek{x}{y},\spvek{x_n'}{u_n}\rangle_{\calX_h\times\calU^*,\calX_h^*\times\calU}.
\end{align*}
But this is a~contradiction to boundedness of the sequence $\big(\spvek{x_n'}{u_n}\big)$ in $\calX_h^*\times \calU$.

{\em Step~3:} We conclude that the desired result holds, i.e., $\overline{\widetilde{ M}}$ as in \eqref{eq:ovtildM} is a~dissipation node on $(\calX_h,\calU)$. Closedness and dissipativity of $\overline{\widetilde{ M}}$ follows from step~1. Closedness of \[F_c\&G_c:\dom \overline{\widetilde{ M}}\supset\calX_h^*\times \calU\to\calX_h\] has been proven in step~2. Further, the assumption \eqref{prop:Mclosb} implies directly that $\overline{\widetilde{ M}}$ has property \eqref{def:dissnode3} in Definition~\ref{def:dissnode}. Finally, we can directly conclude from assumption \eqref{prop:Mclosa}, that \[F_c-\lambda\calR_{\calX_h}:\calX_h^*\supset\dom F_c\to\calX_h\] has dense range. This completes the proof of $\overline{\widetilde{ M}}$ being a~dissipation node.
\end{proof}

\begin{rem}[Dissipation nodes on $(\calX,\calU)$]
The density assumption on $\big(F-\lambda\calR_{\calX_h}^{-1}\big)\calV$ is imposed to guarantee maximal dissipativity of $F_c$, which is -- by Proposition~\ref{prop:maxdiss} -- equivalent to the dissipativity of $F_c^*$. It is an open problem for the authors whether the latter is also implied by density of
\[\calW:=\setdef{x\in \calX_h^*\cap\dom F^*}{F^*x\in \calX_h}\;\text{ in $\calX_h$.}\]
\end{rem}

\section{Applications}\label{sec:appl}

Next we show that several types of partial differential equations with boundary control and observation fit into our framework of port-Hamiltonian systems.

\subsection{An advection-diffusion equation with Dirichlet boundary control}
Let $\Omega\subset\R^d$, $d\in\N$, be a bounded Lipschitz domain with outward normal $\nu:\partial\Omega\to\R^d$ on the boundary, let $a\in L^\infty(\Omega;\C^{d\times d})$ be pointwise Hermitian and positive definite with $a^{-1}\in L^\infty(\Omega;\C^{d\times d})$, and let $b\in L^\infty(\Omega;\R^{d})$ be a~real and divergence-free function with trivial normal trace, i.e., ${\rm div}\,b=0$, and $\nu^\top b=0$ on $\partial\Omega$.\\
Consider the advection-diffusion equation with Dirichlet boundary control and Neumann boundary observation, i.e.,
\begin{equation}\label{eq:advdiff}\begin{aligned}
\tfrac{\partial}{\partial t}x(t,\xi)=&\,\div a(\xi)\,\nabla x(t,\xi)+b(\xi)^\top \nabla x(t,\xi),\quad &&t\ge0, \xi\in\Omega,\\[1mm]
u(t,\xi)=&\,x(t,\xi),\quad y(t,\xi)=\nu(\xi)^\top  a(\xi)\nabla x(t,\xi),&&t\ge0,\, \xi\in\partial\Omega.
\end{aligned}
\end{equation}
We show that this system fits into our framework of port-Hamiltonian systems. Let $\calX=L^2(\Omega)$, and, as Hamiltonian, we take the most standard one
$\mathcal{H}(x)=\frac12\|x\|^2_{L^2(\Omega)}$. That is, the underlying sesquilinear form $h$ is simply the standard inner product in $L^2(\Omega)$, and thus
$\dom h=\calX_h=\calX_h^*=L^2(\Omega)$ and $\widetilde{H}=\Id_{L^2(\Omega)}$.\\ To properly introduce the right formulation and spaces for the above boundary
control problem, consider the {\em trace operator} $\gamma: W^{1,2}(\Omega)\to W^{1/2,2}(\partial\Omega)$ which maps $x\in W^{1,2}(\Omega)$ to its boundary
trace $x|_{\partial\Omega}$, where $W^{1/2,2}(\partial\Omega)$ denotes the Sobolev space of fractional order $1/2$ \cite{adams2003sobolev}. By the {\em trace
theorem} \cite[Thm.~1.5.1.3]{Gris85}, $\gamma$ is bounded and surjective. Further, $H(\div,\Omega)$ is the space of all square integrable functions
whose weak divergence exists and is square integrable. That is, with $W^{1,2}_0(\Omega):=\ker\gamma$,
\[z=\div x\quad\Longleftrightarrow\quad \forall \varphi\in
W^{1,2}_0(\Omega):\;-\langle\nabla\varphi,x\rangle_{L^2(\Omega;\C^d)}=\langle\varphi,z\rangle_{L^2(\Omega)}.\]
Defining $W^{-1/2,2}(\partial\Omega):=W^{1/2,2}(\partial\Omega)^*$, the {\em normal trace} of $x\in H(\div,\Omega)$ is well-defined by $w=\gamma_N x\in
W^{-1/2,2}(\partial\Omega)$ with
\[
 \forall z\in W^{1,2}(\Omega):\;\langle w,\gamma z\rangle_{W^{-1/2,2}(\partial\Omega),W^{1/2,2}(\partial\Omega)}=\langle \div
 x,z\rangle_{L^2(\Omega)}+\langle x,\nabla z\rangle_{L^2(\Omega;\C^d)}.\]
 Green's formula \cite[Chap.~16]{Tart07} yields that, indeed $w(\xi)=\nu(\xi)^\top x(\xi)$ for all $\xi\in\partial\Omega$, if $\partial\Omega$ and $x$ are smooth.
 Further, $\gamma_N:H(\div,\Omega)\to W^{-1/2,2}(\partial\Omega)$ is bounded and surjective \cite[Lem.~20.2]{Tart07}. By boundedness of $\Omega$ and the fact that
 smooth functions are contained in both $L^1(\Omega)$ and $L^2(\Omega)$, we have that $L^2(\Omega)$ is dense in $L^1(\Omega)$.
  As a~consequence, our condition on $b$ is equivalent to
\begin{equation}b\in L^{\infty}(\Omega)\;\;\wedge\;\; \forall\,z\in W^{1,1}(\Omega):\; {\langle\nabla z,b\rangle_{L^1(\Omega;\C^d),L^\infty(\Omega;\C^d)}}=0,\label{eq:divorth}\end{equation}
where the latter stands for the canonical duality product of $L^1$ and $L^\infty\cong (L^1)^*$.
Note that, by boundedness of $\Omega$, we have $L^{\infty}(\Omega)\subset L^{2}(\Omega)$, and thus $b\in H(\div,\Omega)$.\\
Now we introduce the dissipation node corresponding to the advection-diffusion equation with Dirichlet boundary control and Neumann observation, namely
$M=\sbvek{F\&G}{K\&L}$ with
\begin{subequations}\label{eq:advdiff_dissnode}
\begin{multline}
\dom M=\dom{F\&G}\\:=\setdef{\spvek{x}{u}\in W^{1,2}(\Omega)\times W^{1/2,2}(\partial\Omega)}{\,a\nabla x\in
H(\div,\Omega)\;\wedge\;\gamma x=u}
\end{multline}
and
\begin{align}
\forall \, \spvek{x}{u}\in\dom M:\quad F\&G\spvek{x}{u}&=\div\big(a\nabla x\big)+b^\top\nabla x,\;\; K\&L\spvek{x}{u}=-\gamma_N(a\nabla x).
\end{align}
\end{subequations}

\begin{prop}[Dissipation node for advection-diffusion equation with Dirichlet control and Neumann observation]
Let $\Omega\subset\R^d$, $d\in\N$ be a bounded Lipschitz domain. Further, let $a\in L^\infty(\Omega;\C^{d\times d})$ be pointwise Hermitian and positive definite with $a^{-1}\in L^\infty(\Omega;\C^{d\times d})$, and $b\in L^\infty(\Omega;\R^{d})$ with ${\rm div}\,b=0$ and $\gamma_N b=0$. Then $M$ as in \eqref{eq:advdiff_dissnode} is a~dissipation node on $(\calX_h,\calU)$ with $\calX_h=\calX_h^*=L^2(\Omega)$ and $\calU=W^{1/2,2}(\partial\Omega)$, $\calU^*=W^{-1/2,2}(\partial\Omega)$.
\end{prop}
\begin{proof}
We successively show that $M$ fulfills \eqref{def:dissnode1}--\eqref{def:dissnode4} in Definition~\ref{def:dissnode}.
\begin{itemize}
\item[\eqref{def:dissnode1}] {\em Step~1:} We show that $M$ is dissipative. Let $\spvek{x}{u}\in\dom M$.
The product rule for weak derivatives \cite[Thm.~4.25]{Alt16} gives $x\overline{x}\in W^{1,1}(\Omega)$ with
\[\tfrac12\nabla (x\overline{x})=\Re(x \nabla\overline{x})\in L^1(\Omega;\C^d),\]
and thus, by using that $b$ is real-valued,
\begin{align*}
\Re\langle \,x,b^\top\nabla x\rangle_{L^2(\Omega)}
&=\langle \,\Re(x\nabla \overline{x}),b\rangle_{L^1(\Omega;\C^d),L^\infty(\Omega;\C^d)}\\
&=\tfrac12\langle \nabla (x\overline{x}),b\rangle_{_{L^1(\Omega;\C^d),L^\infty(\Omega;\C^d)}}\stackrel{\eqref{eq:divorth}}{=}0.
\end{align*}
Invoking that $\spvek{x}{u}\in\dom M$, we obtain $u=\gamma x$, and thus (for sake of brevity, we leave out the subindices indicating the spaces)
\begin{align*}
\Re\left\langle \spvek{x}{u},M \spvek{x}{u}\right\rangle
&{=}\Re\left\langle \spvek{x}{u},\spvek{\div(a\nabla x)+b^\top\nabla x}{-\gamma_N(a\nabla x)}
\right\rangle\\
&{=}\Re\langle x,\div(a\nabla x)\rangle+\Re\langle x,b^\top\nabla x\rangle
-\langle u,\gamma_N(a\nabla x)\rangle\\
&{=}\Re\langle x,\div(a\nabla x)\rangle+\underbrace{\Re\langle {b}\,x,\nabla x\rangle}_{=0}
-\Re\langle \underbrace{u}_{=\gamma x},\gamma_N(a\nabla x)\rangle\\
&{=}-\Re\langle \nabla x,a\nabla x\rangle+
\Re\langle \gamma_N(a\nabla x),\gamma x\rangle-\Re\langle \gamma x,\gamma_N(a\nabla x)\rangle\\
&{=}-\Re\langle \nabla x,a\nabla x\rangle\leq0,
\end{align*}
where the latter holds by the pointwise positive definiteness of $a$.\\
{\em Step~2:} We show that $ M-\sbmat{\Id}{0}{0}{\calR_\calU}$ is onto. Let $z\in L^2(\Omega)$ and $w\in W^{-1/2}(\partial\Omega)$.
By the Lax-Milgram lemma \cite[Lem.~2.2.1.1]{Gris85}, there exists some $x\in W^{1,2}(\Omega)$, such that
\begin{multline*}
\forall\,\varphi\in W^{1,2}(\Omega):\quad
\langle \nabla\varphi,\nabla x\rangle_{L^2(\Omega;\C^d)}
+ \langle \varphi, x\rangle_{L^2(\Omega)}-\langle \varphi, b^\top \nabla x\rangle_{L^2(\Omega)}+
\langle \gamma\varphi,a\gamma x\rangle_{W^{1/2,2}(\partial\Omega)}\\=-\langle \varphi,z\rangle_{L^2(\Omega)}-\langle \gamma\varphi,w\rangle_{W^{1/2,2}(\partial\Omega),W^{-1/2,2}(\partial\Omega)}.
\end{multline*}
The latter means, by the definition of the weak derivative and the Neumann trace, that, for $u=\gamma x\in W^{1/2,2}(\partial\Omega)$,
\[\big( M-\sbmat{\Id}{0}{0}{\calR_\calU}\big)\spvek{x}{u}=\spvek{z}{w}.\]
{\em Step~3:} We show that $M$ is closed.\\
Since, by steps~1\&2 $M$ is dissipative, and $ M-\sbmat{\Id}{0}{0}{\calR_\calU}$ is onto, the latter operator has a~bounded inverse with operator norm being less or equal to one. This shows that $ M-\sbmat{\Id}{0}{0}{\calR_\calU}$ is closed, and so $ M$ is, as well.
\item[\eqref{def:dissnode2}] We show that $F\&G:\dom M\to L^2(\Omega)$ is closed. Since, by \eqref{def:dissnode1}, $M$ is closed, it suffices, in view of Lemma~\ref{lem:closedblock}, to show that $K\&L\in L(\dom(F\&G),W^{-1/2,2}(\partial\Omega))$. For this, consider sequences $(x_n)$ and $(u_n)$, such that
\begin{align*}
\left(\spvek{x_n}{u_n}\right)&\quad\text{ is bounded in $L^2(\Omega)\times W^{1/2,2}(\partial\Omega)$},\\
\left(F\&G\spvek{x_n}{u_n}\right)&\quad\text{ is bounded in $L^2(\Omega)$.}%\\
%\left(\|K\&L\spvek{x_n}{u_n}\|_{W^{-1/2,2}(\partial\Omega)}\right)&\to \infty.
\end{align*}
We have to show that $(K\& L\spvek{x_n}{u_n})$ is bounded in $W^{-1/2,2}(\partial\Omega)$.

{\em Step~1:} We first prove that $(x_n)$ is bounded in $W^{1,2}(\Omega)$. Since the trace operator $\gamma:W^{1,2}(\Omega)\to W^{1/2,2}(\partial\Omega)$ is bounded and surjective, it possesses a~bounded right inverse $\gamma^-$. That is, $\gamma^-\in L(W^{1/2,2}(\partial\Omega),W^{1,2}(\Omega))$ with $\gamma\gamma^-=\Id_{W^{1/2,2}(\partial\Omega)}$.
Then the sequence $(x_n-\gamma^- u_n)$ is bounded in $L^2(\Omega)$ and, moreover, $x_n-\gamma^-u_n\in W^{1,2}_0(\Omega)$ for all $n\in\N$. By using pointwise positive definiteness and essential boundedness of $a^{-1}$,
\begin{align*}
&\phantom{=}%\langle x-\gamma^- u,z\rangle_{L^2(\Omega)}\\
%&\stackrel{n\to\infty}{\longleftarrow}
\Re\langle x_n-\gamma^-u_n,F\&G \spvek{x_n}{u_n}\rangle_{L^2(\Omega)}\\&=-\Re\langle \nabla x_n-\nabla\gamma^- u_n,a\nabla x_n\rangle_{L^2(\Omega;\C^d)}+\Re\langle x_n-\gamma^- u_n, b^\top\nabla x_n\rangle_{L^{2}(\Omega)}\\
&=-\langle \nabla x_n,a\nabla x_n\rangle_{L^2(\Omega;\C^d)}+\Re\langle \nabla\gamma^- u_n,a\nabla x_n\rangle_{L^2(\Omega;\C^d)}\\&\quad+\Re\langle x_n,b^\top\nabla x_n\rangle_{L^2(\Omega)}-\Re\langle \gamma^-u_n,b^\top\nabla x_n\rangle_{L^2(\Omega)}\\
&\leq-\|a^{-1}\|_{L^\infty(\Omega)}^{-1} \| \nabla x_n\|_{L^2(\Omega;\C^d)}^2+\|\nabla\gamma^- u_n\|_{L^2(\Omega;\C^d)}\|a\|_{L^\infty(\Omega)} \|\nabla x_n\|_{L^2(\Omega;\C^d)}\\&\quad+\|x_n\|_{L^2(\Omega)}\,\|b\|_{L^\infty(\Omega)}\|\nabla x_n\|_{L^2(\Omega;\C^d)}+\| \gamma^-u_n\|_{L^2(\Omega)}\,\|b\|_{L^\infty(\Omega)}\|\nabla x_n\|_{L^2(\Omega;\C^d)}.
\end{align*}
By using that $(x_n)$ is bounded in $L^2(\Omega)$, and $(\gamma^-u_n)$ is bounded in $W^{1,2}(\Omega)$, we obtain the existence of constants  $c_1,c_2>0$, such that, for all $n\in\N$,
\[\| \nabla x_n\|_{L^2(\Omega;\C^d)}^2\leq c_1+c_2\, \| \nabla x_n\|_{L^2(\Omega;\C^d)}.\]
This shows that the real sequence $(\| \nabla x_n\|_{L^2(\Omega;\C^d)})$ is bounded. By further invoking that $(x_n)$ converges in $L^2(\Omega)$, we see that $(x_n)$ is a~bounded sequence in $W^{1,2}(\Omega)$.

\medskip
\noindent{\em Step~2:} We show that $(K\&L\spvek{x_n}{u_n})$ is a~bounded sequence in $W^{-1/2,2}(\partial\Omega)$. To this end, by the
Banach-Steinhaus theorem \cite[Thm.~7.3]{Alt16}, it suffices to show that, for all $y\in W^{ 1/2,2}(\partial\Omega)$, the scalar sequence
\[\left(\left\langle y,K\&L\spvek{x_n}{u_n} \right\rangle_{W^{1/2,2}(\partial\Omega),W^{-1/2,2}(\partial\Omega)}\right)\]%=\left(\left\langle
%y,\gamma_N(a\nabla x_n) \right\rangle_{W^{1/2,2}(\partial\Omega),W^{-1/2,2}(\partial\Omega)}\right)\]
is bounded. Having $y\in W^{1/2,2}(\partial\Omega)$, surjectivity of the trace operator implies that there exists some $w\in W^{1,2}(\Omega)$ with $\gamma w=y$. Then
    \begin{align*}
%\langle w,z\rangle_{L^2(\Omega)}
%&\stackrel{n\to\infty}{\longleftarrow}
\Re\langle w,F\&G \spvek{x_n}{u_n}\rangle_{L^2(\Omega)}&=-\langle \nabla w,a\nabla x_n\rangle_{L^2(\Omega;\C^d)}+\langle w, b^\top\nabla x_n\rangle_{L^{2}(\Omega)}\\&\qquad+\langle y, \gamma_N(a\nabla x_n)\rangle_{W^{1/2,2}(\partial\Omega),W^{-1/2,2}(\partial\Omega)}.
\end{align*}
Now using that, by step~1, $(x_n)$ is a~bounded sequence in $W^{1,2}(\Omega)$, we obtain that the sequence formed by the duality product of $y$ and $\gamma_N(a\nabla x_n)$ has to be bounded. Since $\gamma_N(a\nabla x_n)=K\&L\spvek{x_n}{u_n}$, the result is shown.
%
%the $(\gamma_N x_n)$ is a~bounded sequence in $W^{-1/2,2}(\partial\Omega)$. Assuming the converse
%
%    which shows that $(x_n-v)$ converges to zero in $W^{1,2}(\Omega)$. Further, boundedness of the trace operator yields that the sequence with $u_n-w=\gamma (x_n-z)$ converges to zero in $W^{1/2,2}(\partial\Omega)$. This gives, altogether, $x=v$ and $u=w$. In particular,
%    $\spvek{x}u\in\dom M$ with $(F-\Id)\&G\spvek{x}u=z$, which shows that $(F-\Id)\&G$ is closed.
%    \[=\]
%
%
%     Hence, for all $w\in W^{1,2}(\Omega)$, we have
%    \begin{multline*}
%\langle w,z\rangle_{L^2(\Omega)}\stackrel{n\to\infty}{\longleftarrow}\langle w,(F-\Id)\&G \spvek{x_n}{u_n}\rangle_{L^2(\Omega)}\\
%=-\langle w,x_n\rangle_{L^2(\Omega)}-\langle \nabla w,a\nabla x_n\rangle_{L^2(\Omega;\C^d)}+\langle \gamma w,\gamma x_n\rangle_{L^{2}(\partial\Omega)}.\end{multline*}
%    Now using boundedness of the trace operator $\gamma:W^{1,2}(\Omega)\to W^{1/2,2}(\Omega)$ together with zero convergence of $(x_n)$ in $W^{1,2}(\Omega)$, we obtain that $\langle w,z\rangle_{L^2(\Omega)}=0$ for all $w\in W^{1,2}(\Omega)$, and, by density of $W^{1,2}(\Omega)$ in $L^2(\Omega)$, we have $z=0$.
\item[\eqref{def:dissnode3}] We show that for all $u\in W^{1/2,2}(\Omega)$, there exists some $x\in L^2(\Omega)$ with $\spvek{x}{u}\in \dom M$. Let $u\in W^{1/2,2}(\partial\Omega)$. Since there exists some $x_D\in W^{1,2}(\Omega)$ with $u=\gamma x_D$, the Lax-Milgram lemma further implies that there exists some $x_0\in W^{1,2}_0(\Omega)$, such that
\begin{multline*}
\forall\,\varphi\in W^{1,2}_0(\Omega):\quad
\langle \nabla\varphi,a\nabla x_0\rangle_{L^2(\Omega;\C^d)}+\langle \varphi, x_0\rangle_{L^2(\Omega)}-\langle \varphi, b^\top\nabla x_0\rangle_{L^2(\Omega)}\\=
-\langle \nabla\varphi,a\nabla x_D\rangle_{L^2(\Omega;\C^d)}-\langle \varphi, x_D\rangle_{L^2(\Omega)}+\langle \varphi, b^\top\nabla x_D\rangle_{L^2(\Omega)}.
\end{multline*}
Then $x=x_0+x_D$ fulfills $\gamma x=u$ and
\[
\forall\,\varphi\in W^{1,2}_0(\Omega):\quad
-\langle \nabla\varphi,a\nabla x\rangle_{L^2(\Omega;\C^d)}+\langle \varphi, b^\top\nabla x\rangle_{L^2(\Omega)}=\langle \varphi, x\rangle_{L^2(\Omega)}.
\]
Since the latter means that $\div a\nabla x=x\in L^2(\Omega)$, we obtain that $\spvek{x}u\in\dom F\&G$.
\item[\eqref{def:dissnode4}] We show that $F-\Id$ is surjective. Let $z\in L^2(\Omega)$. Again using the Lax-Milgram lemma, we obtain that there exists some $x\in W^{1,2}_0(\Omega)$ with
\[
\forall\,\varphi\in W^{1,2}_0(\Omega):\quad
\langle \nabla\varphi,\nabla x\rangle_{L^2(\Omega;\C^d)}+\langle \varphi, x\rangle_{L^2(\Omega)}-\langle \varphi, b^\top\nabla x\rangle_{L^2(\Omega)}=
-\langle \varphi,z\rangle_{L^2(\Omega)}.
\]
Since the latter implies that $\div a\nabla x\,- x=z$ together with $\gamma x=0$, the statement is proven.
\end{itemize}
\end{proof}

\begin{rem}[Dirichlet control, Neumann observation, and well-posedness]
Note that the system \eqref{eq:advdiff} is not well-posed in the sense of Definition~\ref{def:wp}, since, in particular, the mapping $\frakB_T$ in Remark~\ref{rem:wp} is unbounded \cite{Schw20}.
\end{rem}

\subsection{The class of {\sc Jacob} and {\sc Zwart}}\label{sec:bjhz}

The systems type addressed in the textbook \cite{Jacob2012} (see also \cite{Augner19,AugnerDis,AugJac14,BastCoro16}) is essentially consisting of a~class of boundary controlled hyperbolic systems in one spatial variable. This type contains a~variety of practically relevant systems, such as the one-dimensional wave equation, the telegraph equations, and the Timoshenko beam.

%Typically these are considered in a bounded interval $[a,b]\subset\R$. We may instead consider $[a,b]=[a,b]$ without loss of generality.\\
Let $H\in L^\infty\left([a,b];\C^{m\times m}\right)$, $m\in\N$, $a,b\in\R$, $a<b$, be pointwise Hermitian and positive definite with $H^{-1}\in L^\infty\left([a,b];\C^{m\times m}\right)$. For $\calX=L^2([a,b];\C^m)$, consider the Hamiltonian $\calH:\calX\to\R$ with
\begin{equation}\label{eq:1dhypHam}
\calH(x)=\tfrac12\,\langle x,H x\rangle_{L^2([a,b];\C^m)}=\tfrac12\,\int_a^b x(\xi)^* H(\xi){x(\xi)}{\rm d}{\xi}.
\end{equation}
Hereby, $H$ is called the {\em Hamiltonian density}.\\
Assume that $P_0\in L^\infty\left([a,b];\C^{m\times m}\right)$ is pointwise dissipative, and let $P_1\in\C^{m\times m}$ be Hermitian and invertible.
Consider the partial differential equation
\begin{align}\label{eq:BH1}
\tfrac{\partial}{\partial t}x(t,\xi)=&\,P_0(\xi)H(\xi)x(t,\xi)+ P_1\tfrac{\partial}{\partial \xi}\big(H(\xi)x(t,\xi)\big),\qquad t\ge0, \xi\in[a,b].
%,\\[1mm]
%u(t,\xi)=&\,x(t,\xi),\\ y(t,\xi)=&\,n(\xi)^\top \nabla a(\xi)x(t,\xi),&&t\ge0,\, \xi\in\partial\Omega.
\end{align}
To construct inputs and outputs, let $W_B,W_C\in\C^{m\times 2m}$ be such that
\begin{equation}\sbmat{0}{\Id_m}{\Id_m}{0}= \sbvek{W_B}{W_C}\sbmat{0}{\Id_m}{\Id_m}{0}\sbvek{W_B}{W_C}^*.\label{eq:PHio}\end{equation}
We consider inputs and outputs formed by suitable boundary evaluations. These are -- in a formal manner -- given by
\begin{equation}\label{eq:BH2}
u(t)=\tfrac1{\sqrt2}W_B\sbmat{P_1}{-P_1}{\Id_m}{\Id_m}\spvek{H(b)\,x(t,b)}{H(a)\,x(t,a)},\quad
y(t)=\tfrac1{\sqrt2}W_C\sbmat{P_1}{-P_1}{\Id_m}{\Id_m}\spvek{H(b)\,x(t,b)}{H(a)\,x(t,a)}.
\end{equation}According to the findings in Section~\ref{sec:quham}, the spaces $\calX_h$, $\calX_h^*$ are equipped with the inner products
\begin{subequations}\label{eq:indualBH}\begin{align}
\langle x,z\rangle_{\calX_h}&=\langle x,H z\rangle_{L^2([a,b];\C^m)}=\int_a^b z(\xi)^* H(\xi){x(\xi)}{\rm d}{\xi},\label{eq:indualBH1}\\
\langle x,z\rangle_{\calX_h^*}&=\langle x,H^{-1} z\rangle_{L^2([a,b];\C^m)}=\int_a^b z(\xi)^* H(\xi)^{-1}{x(\xi)}{\rm d}{\xi},\label{eq:indualBH2}
\end{align}
whereas the duality product is simply the inner product in $L^2([a,b];\C^m)$, i.e.,
\begin{align}
\langle x,z\rangle_{\calX_h^*,\calX_h}&=\langle x,z\rangle_{L^2([a,b];\C^m)}=\int_a^b z(\xi)^*{x(\xi)}{\rm d}{\xi}.\label{eq:indualBH3}
\end{align}
\end{subequations}
By the assumption that $H$ and $H^{-1}$ are bounded, we have $\calX_h=\calX_h^*=\calX=L^2([a,b];\C^m)$. The input and output spaces are given by
$\calU=\calU^*=\C^m$,  and we consider the dissipation node $M=\sbvek{F\&G}{K\&L}$ with
\begin{subequations}\label{eq:BH_dissnode}
\begin{multline}\dom M=\dom{F\&G}\\:=\setdef{\spvek{z}{u}\in W^{1,2}([a,b];\C^m)\times \C^m}{\tfrac1{\sqrt2}W_B\sbmat{P_1}{-P_1}{\Id_m}{\Id_m}\spvek{z(b)}{z(a)}=u}\end{multline}
and
\begin{align}
\forall \, \spvek{z}{u}\in\dom M:\quad F\&G\spvek{z}{u}&=P_0z+ P_1\tfrac{\partial}{\partial \xi}z,\\ K\&L\spvek{z}{u}&=-\tfrac1{\sqrt2}W_C\sbmat{P_1}{-P_1}{\Id_m}{\Id_m}\spvek{z(b)}{z(a)}.
\end{align}
\end{subequations}

\begin{prop}[Dissipation node for hyperbolic systems on a~spatial interval]\label{prop:BHdiss}
Let $a,b\in\R$ with $a<b$, let $m\in\N$, and let $H\in L^\infty([a,b];\C^{m\times m})$ be pointwise Hermitian and positive definite with $H^{-1}\in L^\infty(\Omega;\R^{m\times m})$. Let $P_0\in L^\infty([a,b];\C^{m\times m})$ be pointwise dissipative, and let $P_1\in\C^{m\times m}$ be Hermitian and invertible. Further, let $W_B,W_C\in\C^{m\times 2m}$ such that \eqref{eq:PHio} holds.\\
Then $M$ as in \eqref{eq:BH_dissnode} is a~dissipation node on $(\calX_h,\calU)$ with $\calX_h=\calX_h^*=L^2([a,b];\C^m)$ and $\calU=\calU^*=\C^m$, where $\calX_h$, $\calX_h^*$ are equipped with the inner products and duality product as in \eqref{eq:indualBH}. %, and the inner products in $\calU$ and $\calU^*$ (as well as the duality product in $\calU$) are given by the standard Euclidean inner product in $\C^d$.
\end{prop}
\begin{proof}
We obtain from \eqref{eq:PHio} that ${W_B}\sbmat{0}{\Id_m}{\Id_m}{0}{W_B}^*=0$ and that $W_B$ has full row rank. Then \cite[Thm~2.37]{Vill07} yields that $ M$ is a~system node on $(\C^m,L^2([a,b];\C^m),\C^m)$. Now, by invoking Remark~\ref{rem:dissnode}~\eqref{rem:dissnode4} and Proposition~\ref{prop:dissnode}, it suffices to show that $M:L^2([a,b];\C^m)\supset\dom M\to L^2([a,b];\C^m)$ is dissipative. To this and, let $\spvek{z}{u}\in\dom M$.
Since $P_1$ is constant and Hermitian, integration by parts gives
%\begin{align*}
%\langle z,\tfrac{\partial}{\partial \xi}P_1z\rangle
%&=-\langle \tfrac{\partial}{\partial \xi}z,P_1z\rangle+z(1)^*P_1z(1)-z(0)^*P_1z(0)\\
%&=-\langle \tfrac{\partial}{\partial \xi}P_1z,z\rangle+z(1)^*P_1z(1)-z(0)^*P_1z(0),
%\end{align*}
%which gives
\[\Re\langle z,\tfrac{\partial}{\partial \xi}P_1z\rangle_{L^2([a,b];\C^m)}=\tfrac12 \left.z(\xi)^*P_1z(\xi)\right|_a^b.\]
We can infer from \eqref{eq:PHio} that $\sbvek{W_B}{W_C}\in\C^{2m\times 2m}$ is invertible. Hence, for \[y:=\tfrac1{\sqrt2}W_C\sbmat{P_1}{-P_1}{\Id_m}{\Id_m}\spvek{z(b)}{z(a)},\] and by invoking that \[u=\tfrac1{\sqrt2}W_B\sbmat{P_1}{-P_1}{\Id_m}{\Id_m}\spvek{z(b)}{z(a)},\]
we obtain
\begin{align*}
\Re\langle \tfrac{\partial}{\partial \xi}z,P_1z\rangle_{L^2([a,b];\C^m)}
&=\tfrac12 \left.z(\xi)^*P_1z(\xi)\right|_a^b\\
&=\tfrac14\left( \sbmat{P_1}{-P_1}{\Id_m}{\Id_m}\spvek{z(a)}{z(b)}\right)^*\left(\sbmat{0}{\Id_m}{\Id_m}{0}\sbmat{P_1}{-P_1}{\Id_m}{\Id_m}\spvek{z(a)}{z(b)}\right)\\
&=\tfrac14\left( \sqrt2\sbvek{W_B}{W_C}^{-1}\spvek{u}{y}\right)^*\left(\sbmat{0}{\Id_m}{\Id_m}{0}\sqrt2\sbvek{W_B}{W_C}^{-1}\spvek{u}{y}\right)\\
%&=\langle \spvek{u}{y},\sbvek{W_B}{W_C}^{-*}\sbmat{0}{\Id_d}{\Id_d}{0}\sbvek{W_B}{W_C}^{-1}\spvek{u}{y}\rangle\\
&=\tfrac12\spvek{u}{y}^*\left(\sbvek{W_B}{W_C}\sbmat{0}{\Id_m}{\Id_m}{0}\sbvek{W_B}{W_C}^{*}\right)^{-1}\spvek{u}{y}\\
&=\tfrac12 \spvek{u}{y}^*\sbmat{0}{\Id_m}{\Id_m}{0}^{-1}\spvek{u}{y}=\tfrac12 \spvek{u}{y}^*\sbmat{0}{\Id_m}{\Id_m}{0}\spvek{u}{y}=\Re({u}^*{y}).
\end{align*}
Then, dissipativity follows from
\begin{align*}
&\phantom{=}\Re\left\langle \spvek{z}{u},M \spvek{x}{u}\right\rangle_{L^2([a,b];\C^m)\times \C^m}\\
&{=}\Re\left\langle \spvek{z\\[-1mm]}{u},\spvek{P_0z+ P_1\tfrac{\partial}{\partial \xi}z}{-y}
\right\rangle_{L^2([a,b];\C^m)\times \C^m}\\
&{=}\Re\langle z,P_0z\rangle_{L^2([a,b];\C^m)}+\Re\langle z,\tfrac{\partial}{\partial \xi}P_1z\rangle_{L^2([a,b];\C^m)}-\Re( u^*y)\\
&{=}\Re\langle z,P_0z\rangle_{L^2([a,b];\C^m)}\leq0,
\end{align*}
where the latter holds due to pointwise dissipativity of $P_0$.
\end{proof}

\begin{rem}[Hyperbolic systems on a~spatial interval]\
\begin{enumerate}[(a)]
\item It is shown in \cite[Sec.~13.2]{Jacob2012} that, under the assumptions made throughout this section, \eqref{eq:BH1}\&\eqref{eq:BH2} forms a well-posed system in the sense of Definition~\ref{def:wp}.
\item In \cite{Jacob2012}, some quite more general input-output boundary configurations are considered, which lead to an energy balance of the form
\[
\calH(x(t))-\calH(x_0)\leq\int_0^t \spvek{u(\tau)}{y(\tau)}^*S\spvek{u(\tau)}{y(\tau)}{\rm d}\tau\]
for some Hermitian matrix $S\in\C^{2m\times 2m}$. Note that such ports (i.e., the collection of input and output) are not port-Hamiltonian in the classical sense of \cite{van2014port}.
\end{enumerate}
\end{rem}

For later use we record the following lemma on the adjoint of $ M$ as in \eqref{eq:BH_dissnode}, which is -- loosely speaking -- obtained from $ M$ by reflecting the spatial interval $[a,b]$. The proof is analogous to \cite[Prop.~3.4.3]{AugnerDis}, where higher order systems with $u=0$ are considered. It is therefore omitted.

\begin{lem}[Adjoint of the dissipation node \eqref{eq:BH_dissnode}]\label{lem:BHadj}
Let $a,b\in\R$ with $a<b$, let $m\in\N$, and let $H\in L^\infty([a,b];\C^{m\times m})$ be pointwise Hermitian and positive definite with $H^{-1}\in L^\infty(\Omega;\R^{m\times m})$. Let $P_0\in L^\infty([a,b];\C^{m\times m})$ be pointwise dissipative, and let $P_1\in\C^{m\times m}$ be Hermitian and invertible. Further, let $W_B,W_C\in\C^{m\times 2m}$ such that \eqref{eq:PHio} holds.\\
Then, for $\calX=\calX^*=L^2([a,b];\C^{m})$, the adjoint of $ M:\calX\times\calU\to\calX\times\calU^*$ as in \eqref{eq:BH_dissnode} is given by $M^*=\left[\begin{smallmatrix}{[F\& G]^d}\\{[K\& L]^d}\end{smallmatrix}\right]$ with
\begin{multline*}\dom M^*=\dom{[F\& G]^d}\\:=\setdef{\spvek{z}{v}\in W^{1,2}([a,b];\C^m)\times \C^m}{\tfrac1{\sqrt2}W_B\sbmat{P_1}{-P_1}{\Id_m}{\Id_m}\spvek{z(a)}{z(b)}=v}\end{multline*}
and
\begin{align*}
\forall \, \spvek{z}{v}\in\dom M^*:\quad [F\& G]^d\spvek{z}{v}&=P_0^*z- P_1\tfrac{\partial}{\partial \xi}z,\\ [K\& L]^d\spvek{z}{v}&=-\tfrac1{\sqrt2}W_C\sbmat{P_1}{-P_1}{\Id_m}{\Id_m}\spvek{z(a)}{z(b)}.
% M&=\sbvek{F\&G}{K\&L}.
\end{align*}
\end{lem}

\begin{rem}[Partially homogeneous boundary conditions]\label{rem:parthom}
In \cite{Jacob2012}, also input-output configurations together with some additional homogeneous boundary conditions have been considered. That is, for some $m_1\leq m$, $W_{B,1},W_{C,1}\in\C^{m_1\times 2m}$, $W_{B,2}\in\C^{m_2\times 2m}$, $m_2:=m-m_1$, such that $\sbvek{W_{B,1}}{W_{B,2}}$ has full row rank and
\[
\left[\begin{smallmatrix}{0}&0&{\Id_{m_1}}\\0&0&0\\{\Id_{m_1}}&{0}&{0}\end{smallmatrix}\right]= \left[\begin{smallmatrix}W_{B,1}\\W_{B,2}\\W_{C,1}\end{smallmatrix}\right]\bmat{0}{\Id_m}{\Id_m}{0}\left[\begin{smallmatrix}W_{B,1}\\W_{B,2}\\W_{C,1}\end{smallmatrix}\right]^*,
\]
we consider \eqref{eq:BH1} with input and output
\[\begin{aligned}
u(t)&=\tfrac1{\sqrt2}W_{B,1}\sbmat{P_1}{-P_1}{\Id_m}{\Id_m}\spvek{H(b)\,x(t,b)}{H(a)\,x(t,a)},\\
y(t)&=\tfrac1{\sqrt2}W_{C,1}\sbmat{P_1}{-P_1}{\Id_m}{\Id_m}\spvek{H(b)\,x(t,b)}{H(a)\,x(t,a)}
\end{aligned}\]
together with the homogeneous boundary condition
\[0=W_{B,2}\sbmat{P_1}{-P_1}{\Id_m}{\Id_m}\spvek{H(b)\,x(t,b)}{H(a)\,x(t,a)}.\]
Simple linear algebra shows that we can extend by some $W_{C,2}\in\C^{m_2\times 2m}$, such that \eqref{eq:PHio} is fulfilled for $W_B:=\sbvek{W_{B,1}}{W_{B,2}}$, $W_C:=\sbvek{W_{C,1}}{W_{C,2}}$. Clearly, the Hamiltonian for this system is again given by \eqref{eq:1dhypHam}. To formulate the corresponding dissipation node, we first take the one in \eqref{eq:BH_dissnode} (which has yet too many inputs and outputs), and then reduce these by applying Proposition~\ref{prop:iored} with $\calU_{\rm new}=\calU_{\rm new}^*=\C^{m_1}$ and
\[U=\sbvek{\Id_{m_1}}{0}\in\C^{m\times m_1}.\]
\end{rem}

\begin{rem}[Hyperbolic systems of higher order]
In \cite{AugnerDis,Vill07}, a more general class of boundary-controlled hyperbolic class of systems (which for instance allows to incorporate the Euler-Bernoulli beam) in one spatial variable has been considered, namely, for a~Hamiltonian density ${H}\in L^\infty([a,b];\C^{d\times d})$ and $P_0:[a,b]\to\C^{d\times d}$ with assumptions as above, $N\in\N$, $P_k=(-1)^{k+1}P_k^\ast\in\C^{d\times d}$ for $k=1,\ldots,N$ with $P_N$ invertible, consider the partial differential equation
    \begin{align}\label{eq:BH1gen}
\tfrac{\partial}{\partial t}x(t,\xi)=&\,P_0(\xi)H(\xi)x(t,\xi)+ \sum_{k=1}^NP_k\tfrac{\partial^k}{\partial\zeta^k}\big(H(\xi)x(t,\xi)\big),\quad &&t\ge0, \xi\in[a,b].
%,\\[1mm]
%u(t,\xi)=&\,x(t,\xi),\\ y(t,\xi)=&\,n(\xi)^\top \nabla a(\xi)x(t,\xi),&&t\ge0,\, \xi\in\partial\Omega.
\end{align}
For $m:={Nd}$ and $W_B,W_C\in\C^{m\times 2m}$ with
\[\sbmat{0}{\Id_{m}}{\Id_{m}}{0}= \sbvek{W_B}{W_C}\sbmat{0}{\Id_{m}}{\Id_{m}}{0}\sbvek{W_B}{W_C}^*,\]
inputs and outputs are given by
\[u(t)=\tfrac1{\sqrt2}W_B\sbmat{\Lambda}{-\Lambda}{\Id_m}{\Id_m}\spvek{w(t,1)}{w(t,0)},\quad
y(t)=\tfrac1{\sqrt2}W_C\sbmat{\Lambda}{-\Lambda}{\Id_m}{\Id_m}\spvek{w(t,1)}{w(t,0)},
\]
where
		\[w(t,\xi)=\left(\begin{smallmatrix}
		H(\xi)x(t,\xi)\\
		\tfrac{\partial}{\partial\xi}H(\xi)x(t,\xi)\\
		\vdots\\
		\tfrac{\partial^{N-1}}{\partial\xi^{N-1}}H(\xi)x(t,\xi)
		\end{smallmatrix}\right),\quad
\Lambda:=\begin{bmatrix}
		P_1 & P_2 & \cdots & \cdots & P_N\\
		-P_2 & -P_3 & \cdots & -P_N & 0\\
		\vdots & \vdots &\vdots & \vdots & \vdots\\
		(-1)^{N-1}P_N & 0 &\cdots& 0 & 0
		\end{bmatrix}.
		\]
It can indeed be shown that this forms as well a~port-Hamiltonian system in the sense of this article. Due to the high level of technicality and large amount of subindices, this is left out here.
\end{rem}

\subsection{The class of {\sc Jacob} and {\sc Zwart} - extended to singular Hamiltonian densities}\label{sec:bjhz_sing}

We revisit the class treated in the antecedent section. However, we now relax the assumption on the Hamiltonian, such that our presented theory on Hamiltonians in conjunction with quasi Gelfand triples comes truly into operation. Instead of imposing that the Hamiltonian density and its pointwise inverse are bounded, we now assume that, for $m\in\N$, $a,b\in\R$ with $a<b$, the Hermitian positive definite-valued function $H:[a,b]\to\C^{m\times m}$ fulfills
\[H,H^{-1}\in L^1\left([a,b];\C^{m\times m}\right).\]
Loosely speaking, $H$ may now have zeros and singularities of low order. The Hamiltonian is again defined by \eqref{eq:1dhypHam}, which means that the underlying sesquilinear form $h$ is
\begin{align*}
h(x,z)&=\langle x,H z\rangle_{L^2([a,b];\C^m)}=\int_a^b z(\xi)^* H(\xi){x(\xi)}{\rm d}{\xi},\\
\dom h&=\setdef{x\in L^2([a,b];\C^m)}{H^{1/2}x\in L^2([a,b];\C^m)},
\end{align*}
where $H^{1/2}:[a,b]\to\C^{m\times m}$ is the pointwise matrix square root of $H$. Positivity of $h$ follows immediately from the pointwise positive definiteness of $H$. Moreover, it can be seen that
\begin{align*}
 \iota:\quad L^2([a,b];\C^m)&\to\dom h,\\
x&\mapsto (\Id_m+H)^{-1/2}x,
\end{align*}
is an isometric isomorphism. Hence, the form $h$ is closed. Density of $\dom h$ in $L^2([a,b];\C^m)$ follows from \[L^\infty([a,b];\C^m)\subset\dom h,\]
whereat the latter is can be obtained by using $H\in L^1([a,b];\C^{m\times m})$ together with the H\"older inequality.
Since the inner products in $\calX_h$ and $\calX_h^*$ are respectively given by \eqref{eq:indualBH1} and \eqref{eq:indualBH2}, it can be concluded that
\begin{equation}\label{eq:XhBH}
\calX_h = H^{-1/2}\cdot L^2([a,b];\C^m),\quad \calX_h^* = H^{1/2}\cdot L^2([a,b];\C^m).
%\calX_h^*=\setdef{x:[a,b]\to\C^m}{H^{-1/2}x\in L^2([a,b];\C^m)}
%\end{aligned}
\end{equation}
%modulo the subspace of functions vanishing almost everywhere.
By construction of quasi Gel\-fand triples, $\langle x,z\rangle_{\calX_h^*,\calX_h}$ is -- as in \eqref{eq:indualBH3} -- the integral of $z(\xi)^*x(\xi)$ over $[a,b]$. Further, the Riesz isomorphism is given by the multiplication operator
\begin{equation}\begin{aligned}
\calR_{\calX_h}:\qquad \calX_h&\to\calX_h^*,\\
x&\mapsto Hx.
\end{aligned}\label{eq:RieszBH}\end{equation}
A further property is shown in the following auxiliary result.
\begin{lem}[$\calX_h$ and $\calX_h^*$ are contained in $L^1$]\label{eq:XhL1}
Let $a,b\in\R$ with $a<b$, let $m\in\N$, and let $H\in L^1([a,b];\C^{m\times m})$ be pointwise Hermitian and positive definite with $H^{-1}\in L^1([a,b];\C^{m\times m})$. Then the spaces $\calX_h$, $\calX_h^*$ are both contained in $L^1([a,b];\C^{m\times m})$.
\end{lem}
\begin{proof}
By using that, by elementary linear algebra, the trace of $H$ equals to the sum of squares of the entries of $H^{1/2}$, we obtain from $H\in L^1([a,b];\C^{m\times m})$ that $H^{1/2}\in L^2([a,b];\C^{m\times m})$. Then \eqref{eq:XhBH} together with H\"older's inequality gives
\[
z=H^{1/2}\cdot H^{-1/2}z\in L^1([a,b];\C^m)\quad \forall z\in\calX_h^*.
\]
The proof of $\calX_h\subset L^1([a,b];\C^{m\times m})$ is completely analogous, since it follows by the above argumentation in which $H$ is replaced with $H^{-1}$.
\end{proof}
As in Section~\ref{sec:bjhz}, we assume that
 $P_0\in L^\infty\left([a,b];\C^{m\times m}\right)$ is pointwise dissipative, and $P_1\in\C^{m\times m}$ is Hermitian and invertible.
Under our modified assumptions on the Hamiltonian, we consider the spatially one-dimensional hyperbolic partial differential equation \eqref{eq:BH1}. Inputs and outputs are chosen as in the
previous section, as well. That is, for $W_B,W_C\in\C^{m\times 2m}$ with \eqref{eq:PHio}, inputs and outputs are formed by boundary values as in
\eqref{eq:BH2}. In particular, $\calU=\calU^*=\C^m$.

To construct the corresponding dissipation node, we proceed as presented at the end of Section~\ref{Ssec:pH}: We first take the dissipation node $ M$ as in \eqref{eq:BH_dissnode}, which is thereafter restricted to the space $\calZ$ as in \eqref{eq:Zdef}. From this operator we take the closure to end up with an operator
\[\overline{\widetilde{M}}:\calX_h^*\times\calU\supset\dom\overline{\widetilde{M}}\to \calX_h^*\times\calU,\]
which is now shown to be a~dissipation node on $(\calX_h,\calU)$. In particular, we show that this operator obtained by closure is the one which truly takes classical derivatives, if smooth functions are plugged in. First, we denote the space of continuously differentiable functions vanishing at $a$ and $b$ by $C^1_0([a,b];\C^m)$.

\begin{prop}[The dissipation node]\label{prop:BHsing}
Let $a,b\in\R$ with $a<b$, let $m\in\N$, and let $H\in L^1([a,b];\C^{m\times m})$ be pointwise Hermitian and positive definite with $H^{-1}\in L^1([a,b];\C^{m\times m})$. Let $P_0\in L^\infty([a,b];\C^{m\times m})$ be pointwise dissipative, and let $P_1\in\C^{m\times m}$ be Hermitian and invertible. Further, let $W_B,W_C\in\C^{m\times 2m}$ such that \eqref{eq:PHio} is fulfilled. Then the following holds for the operator \[ M=\sbvek{F\&G}{K\&L}:L^2([a,b];\C^m)\times\C^m\supset\dom{ M}\to L^2([a,b];\C^m)\times\C^m\]  as in \eqref{eq:BH_dissnode}, and spaces $\calX_h$, $\calX_h^*$ as in \eqref{eq:XhBH}:
\begin{enumerate}[(a)]
\item\label{prop:BHsing1} For all $u\in\C^m$, there exists some $z\in C^1([a,b];\C^m)$ with $\spvek{z}{u}\in\dom M$.
\item\label{prop:BHsing2} For all $y\in\C^m$, there exists some $w\in C^1([a,b];\C^m)$ with $\spvek{w}{y}\in\dom M^*$.
\item\label{prop:BHsing3} The space $C^1_0([a,b];\C^m)\subset L^2([a,b];\C^m)$ is dense in $\calX_h^*$. Moreover, for
\[
\calW:=\setdef{z\in\calX_h^*\cap C^1([a,b];\C^m)}{\,W_B\sbmat{P_1}{-P_1}{\Id_m}{\Id_m}\spvek{z(b)}{z(a)}=0},
\]
the space
\[(F-\calR_{\calX_h}^{-1})\,\calW\]
is dense in $\calX_h$.
\item The restriction of $ M$ to
\[
\setdef{\spvek{z}{u}\in (\calX_h^*\times\C^m)\cap\dom M}{\,M\spvek{z}{u}\in\calX_h^*}
\]
is closable as an operator defined on a~subspace of $\calX_h^*\times\C^m$, and mapping to $\calX_h\times\C^m$. Its closure
\[
\overline{\widetilde{ M}}=:\sbvek{F_c\&G_c}{K_c\&L_c}:\calX_h^*\times\C^m\supset \dom\overline{\widetilde{ M}}\to\calX_h\times\C^m
\]
is a~dissipation node on $(\calX_h,\C^m)$.
\end{enumerate}
\end{prop}
\begin{proof}\
\begin{enumerate}[(a)]
\item Let $u\in\C^m$. Then the result follows, since, by a~simple interpolation, there exists some $z\in C^1([a,b];\C^m)$ with \[\tfrac1{\sqrt{2}}W_B\sbmat{P_1}{-P_1}{\Id_m}{\Id_m}\spvek{z(b)}{z(a)}=u.\]
\item This follows from the argumentation and the fact that, by Lemma~\ref{lem:BHadj}, the adjoint of $ M$ satisfies the same assumptions as $ M$.
\item The H\"older inequality together with $H^{-1}\in L^1([a,b];\C^{m\times m})$ yields
 $L^\infty([a,b];\C^m)\subset \calX_h^*$, which indeed leads to
$C^1_0([a,b];\C^m)\subset \calX_h^*$. Next we show that this inclusion is dense. To this end, by the Hahn-Banach theorem \cite[Thm.~6.15]{Alt16}, it suffices to show that any $z\in\calX_h$ with
\begin{equation}\langle\varphi ,z\rangle_{\calX_h^*,\calX_h}=0\quad\forall \varphi\in \calW\label{eq:Xhorth}\end{equation}
must be zero.
Assume that $z\in\calX_h$ fulfills \eqref{eq:Xhorth}. By using Lemma~\ref{eq:XhL1}, we have $z\in L^1([a,b];\C^m)$ with
\[\langle\varphi ,z\rangle_{L^\infty([a,b];\C^m),L^1([a,b];\C^m)}=0\quad\forall \varphi\in C^1_0([a,b];\C^m)\subset \calW.\]
The fundamental lemma of calculus of variations \cite[Thm.~6.3-2]{Cia13} now leads to $z=0$.\\
In the sequel, we show that $(F-\calR_{\calX_h}^{-1})\,\calW$ is dense in $\calX_h$. We show the equivalent statement that only $z=0$ fulfills
$0=\langle z,(F-\calR_{\calX_h}^{-1})\varphi\rangle_{\calX_h^*,\calX_h}$ for all $\varphi\in\calW$.\\
Assume that $z\in \calX_h^*$ with
\[0=\langle z,(F-\calR_{\calX_h}^{-1})\varphi\rangle_{\calX_h^*,\calX_h}\quad  \forall \,\varphi\in\calW.\]
Since $\calW\subset L^\infty([a,b];\C^m)$ and, by combining \eqref{eq:RieszBH} with Lemma~\ref{eq:XhL1},
\[H^{-1}z\in\calX_h\subset L^1([a,b];\C^m),\] we have
\begin{align*}
0&=\langle z,(F-\calR_{\calX_h}^{-1})\varphi\rangle_{\calX_h^*,\calX_h}\\
 &=\langle z,P_0\varphi+ P_1\tfrac{\partial}{\partial \xi}\varphi\rangle_{\calX_h^*,\calX_h}-\langle z,\varphi\rangle_{\calX_h^*}\\
 &=\langle P_0^*z,\varphi\rangle_{\calX_h^*,\calX_h}+\langle P_1z,\tfrac{\partial}{\partial \xi}\varphi\rangle_{\calX_h^*,\calX_h}-\langle z,\varphi\rangle_{\calX_h^*}\\
 &=\langle P_0^*z,\varphi\rangle_{L^1([a,b];\C^m),L^\infty([a,b];\C^m)}+\langle P_1z,\tfrac{\partial}{\partial \xi}\varphi\rangle_{L^1([a,b];\C^m),L^\infty([a,b];\C^m)}\\&\quad-\langle H^{-1}z,\varphi\rangle_{L^1([a,b];\C^m),L^\infty([a,b];\C^m)}.
\end{align*}
Since the above equality in particular holds for all $\varphi\in C^1_0([a,b];\C^m)$, and $z,H^{-1}z\in L^1([a,b];\C^m)$, the definition of the weak derivative gives rise to $z\in W^{1,1}([a,b];\C^m)$ with
\[0=P_0^*z-P_1\tfrac{\partial}{\partial \xi}z-H^{-1}z.\]
This implies that $z\in L^\infty([a,b];\C^m)\subset\calX_h$, such that, further, $\tfrac{\partial}{\partial \xi}z\in\calX_h$.
Now an integration by parts yields that for all $\varphi\in\calW$,
%by invoking
%\[W_B\sbmat{P_1}{-P_1}{\Id_m}{\Id_m}\spvek{\varphi(b)}{\varphi(a)}=0,\]
that
\begin{align*}
0&=\langle z,(F-\calR_{\calX_h}^{-1})\varphi\rangle_{\calX_h^*,\calX_h}\\
&=\langle P_0^*z,\varphi\rangle_{L^1([a,b];\C^m),L^\infty([a,b];\C^m)}
+\langle P_1z,\tfrac{\partial}{\partial \xi}\varphi\rangle_{L^1([a,b];\C^m),L^\infty([a,b];\C^m)}\\&\quad-\langle H^{-1}z,\varphi\rangle_{L^2([a,b];\C^m)}\\
&=\langle P_0^*z-P_1\tfrac{\partial}{\partial \xi}z-H^{-1}z,\varphi\rangle_{L^1([a,b];\C^m),L^\infty([a,b];\C^m)}+\left.\varphi(\xi)^*P_1z(\xi)\right|_a^b\\
&=\left.\varphi(\xi)^* P_1z(\xi)\right|_a^b.
\end{align*}
By invoking that $\sbvek{W_B}{W_C}$ is invertible (which is a~simple consequence of \eqref{eq:PHio}) and
    \begin{equation}W_B\sbmat{P_1}{-P_1}{\Id_m}{\Id_m}\spvek{\varphi(b)}{\varphi(a)}=0,\label{eq:phibnd}\end{equation}
we have, for
\begin{align*}
u_z&:=\tfrac1{\sqrt{2}}W_B\sbmat{P_1}{-P_1}{\Id_m}{\Id_m}\spvek{z(b)}{z(a)},\\
y_z&:=\tfrac1{\sqrt{2}}W_C\sbmat{P_1}{-P_1}{\Id_m}{\Id_m}\spvek{z(b)}{z(a)},\qquad
y_\varphi:=\tfrac1{\sqrt{2}}W_C\sbmat{P_1}{-P_1}{\Id_m}{\Id_m}\spvek{\varphi(b)}{\varphi(a)},
\end{align*}
that
\begin{equation}
\begin{aligned}
0&=\left.\varphi(\xi)^* P_1z(\xi)\right|_a^b\\
&=\tfrac12\left(\sbmat0{\Id_m}{\Id_m}0 \sbmat{P_1}{-P_1}{\Id_m}{\Id_m}\spvek{\varphi(b)}{\varphi(a)}\right)^*\left(\sbmat{P_1}{-P_1}{\Id_m}{\Id_m}\spvek{z(b)}{z(a)}\right)\\
&=\left(\sbvek{W_B}{W_C}^{-1}\spvek{0}{y_\varphi}\right)^*\left(\sbmat0{\Id_m}{\Id_m}0 \sbvek{W_B}{W_C}^{-1}\spvek{u_z}{y_z}\right)\\
&\stackrel{\eqref{eq:PHio}}{=}\spvek{0}{y_\varphi}^*\left(\sbmat0{\Id_m}{\Id_m}0 \spvek{u_z}{y_z}\right)\\
&= y_\varphi^*u_z.
\end{aligned}\label{eq:disscalc}
\end{equation}
Invertibility of $\sbmat{P_1}{-P_1}{\Id_m}{\Id_m}$ yields that we can find - by a~simple interpolation argument as in \eqref{prop:BHsing1}\&\eqref{prop:BHsing2} - for all $y_\varphi\in\C^m$ some $\varphi\in C^1([a,b];\C^m)$ with
\[\spvek{\varphi(b)}{\varphi(a)}=\sqrt2 \sbmat{P_1}{-P_1}{\Id_m}{\Id_m}^{-1}\sbvek{W_B}{W_C}^{-1}\spvek{0}{y_\varphi}.\]
Then, by the fact that \eqref{eq:disscalc} holds for all $\varphi\in C^1([a,b];\C^m)$ satisfying \eqref{eq:phibnd}, we are led to $u_z=0$. In other words,  $z\in\dom F$, and by an argumentation as in the proof of Proposition~\ref{prop:BHdiss}, we obtain
\[
0=\Re\langle z,(F-\calR_{\calX_h}^{-1})z\rangle_{\calX_h^*,\calX_h}=\Re\langle (P_0^*-H^{-1})z,z\rangle_{L^1([a,b];\C^m),L^\infty([a,b];\C^m)}.\]
Now pointwise dissipativity of $P_0$ together with pointwise positivity of $H$ leads to $z=0$.
\item We have
\[C^1_0([a,b];\C^m)\subset \calV:=\setdef{x\in\calX_h^*\cap\dom F}{Fx\in \calX_h},\]
whence, by \eqref{prop:BHsing3}, $\calV$ is dense in $\calX^*$. Further, by $\calW\subset\calV$ together with $(F-\calR_{\calX_h}^{-1})\,\calW$ being dense in $\calX_h$, we have that $(F-\calR_{\calX_h}^{-1})\,\calV$ is dense in $\calX_h$, as well. That is, assumption \eqref{prop:Mclosa} in
Proposition~\ref{prop:Mclos} is fulfilled. Since, moreover, \eqref{prop:BHsing1} and \eqref{prop:BHsing2} imply that assumptions \eqref{prop:Mclosb}\&\eqref{prop:Mclosc} in
Proposition~\ref{prop:Mclos} are true for $M$, we have that $\overline{\widetilde{ M}}$ is a~dissipation node on $(\calX_h,\C^m)$.
\end{enumerate}
\end{proof}

\begin{rem}[Partially homogeneous boundary conditions]\label{rem:parthomnew}
The statements of Remark~\ref{rem:parthom} also fully apply to the setup presented in the current section.
\end{rem}

\begin{ex}[Vibrating string]
A possible application is given by a~linear wave equation on a~spatial interval $[a,b]$, $a<b$ \cite{Jacob2012}. In compact notation, this reads
\[\tfrac{\partial}{\partial t}\begin{pmatrix}\mathbf{q}(t,\xi)\\[1mm]\mathbf{p}(t,\xi)\end{pmatrix}=\begin{bmatrix}0&\tfrac{\partial}{\partial \xi}\\\tfrac{\partial}{\partial \xi}&-d(\xi)\end{bmatrix}\begin{pmatrix}T(\xi)\;\;\mathbf{q}(t,\xi)\\[1mm]\rho(\xi)^{-1}\mathbf{p}(t,\xi)\end{pmatrix},\quad t\geq0,\,\xi\in[a,b],
\]
where $T:[a,b]\to\R$ is {Young's modulus}, $\rho:[a,b]\to\R$ is the {\em longitudinal mass density}, and $d:[a,b]\to\R$ stands for a~distributed damping. The time dependent variables $\mathbf{q}$ and $\mathbf{p}$ respectively stand for the strain and momentum.\\
We assume that the string is fixed at the left end of the string, whereas control is done via force at the right end; the output consists of the co-located velocity. That is
\[u(t)=-T(b)\,\mathbf{q}(t,b),\quad y(t)={\rho(b)}^{-1}\mathbf{p}(t,b),\quad 0={\rho(a)}^{-1}\mathbf{p}(t,a).
\]
%
%
%
%\[W_B
%=\tfrac1{\sqrt{2}}\left[\begin{smallmatrix}0&-1&1&0\\0&1&1&0\end{smallmatrix}\right],\quad W_C=\tfrac1{\sqrt{2}}\left[\begin{smallmatrix}\\1&0&0&1\end{smallmatrix}\right],
%\]
%
%
%\[u(t)=
%\left(\begin{smallmatrix}\phantom{-}T(a)\,\mathbf{q}(t,a)\\-T(b)\,\mathbf{q}(t,b)\end{smallmatrix}\right),\quad
%y(t)=\spvek{{\rho(a)}^{-1}\mathbf{p}(t,a)}{}.
%\]
%
Our assumptions are that $\rho$ and $T$ are positive-valued with $\rho,\rho^{-1},T,T^{-1}\in L^1([a,b])$. We further assume that $d\in L^\infty([a,b])$ is nonnegative-valued.\\
This fits in the class treated in this section (in particular, see Remarks~\ref{rem:parthom}\&\ref{rem:parthomnew}), within whose notation we have
\begin{align*}
H&=\sbmat{T}00{\rho^{-1}},\quad P_0=\sbmat000{-d},\quad P_1=\sbmat0110,\\
W_{B,1}&
=\tfrac1{\sqrt{2}}\left[\begin{smallmatrix}0&1&1&0\end{smallmatrix}\right],\quad W_{B,2}=\tfrac1{\sqrt{2}}\left[\begin{smallmatrix}-1&0&0&1\end{smallmatrix}\right],\quad W_{C,1}=\tfrac1{\sqrt{2}}\left[\begin{smallmatrix}1&0&0&1\end{smallmatrix}\right],
\end{align*}
and the energetic space and its dual read
\begin{align*}
\calX_h&=\setdef{\spvek{\mathbf{q}}{\mathbf{p}}}{T^{1/2}\mathbf{q},
\rho^{-1/2}\mathbf{p}\in L^2([a,b])},\quad
\calX_h^*=\setdef{\spvek{\mathbf{F}}{\mathbf{v}}}{T^{-1/2}\mathbf{F},
\rho^{1/2}\mathbf{v}\in L^2([a,b])}.
\end{align*}
In physical terms, the Hamiltonian consists of the sum of potential and kinetic energy, whereas the elements of $\calX_h^*$ are consisting of pairs of forces and velocities. Therefore, $\calX_h$ truly consists of all pairs of strains and momenta which correspond to finite potential and kinetic energy, resp. Note that the assumptions $\rho,T^{-1}\in L^1([a,b])$ mean that the string has finite mass and positive spring rate, resp. The authors are not aware of a~physical interpretation of $\rho^{-1},T\in L^1([a,b])$.\\
By combining Proposition~\ref{prop:BHsing} with Theorem~\ref{thm:pHnode} and Proposition~\ref{prop:solex}, we further obtain insight into the solution behavior of the system. In particular, the operator
\begin{align*}
&A:\calX_h\supset\dom A\to\calX_h, \quad A\spvek{\mathbf{q}}{\mathbf{p}}=\left[\begin{smallmatrix}0&\tfrac{\partial}{\partial \xi}\\\tfrac{\partial}{\partial \xi}&-d(\xi)\end{smallmatrix}\right]\spvek{T\mathbf{q}}{\rho^{-1}\mathbf{p}}\\
&\dom A=\setdef{\spvek{\mathbf{q}}{\mathbf{p}}\in\calX_h}{\left[\begin{smallmatrix}0&\tfrac{\partial}{\partial \xi}\\\tfrac{\partial}{\partial \xi}&-d\end{smallmatrix}\right]\spvek{T\mathbf{q}}{\rho^{-1}\mathbf{p}}\in\calX_h\text{ with }(\rho^{-1}\mathbf{p})(a)=(T\mathbf{q})(b)=0}
\end{align*}
generates a~contractive semigroup on $\calX_h$.
\end{ex}

\begin{rem}[Transmission line and Timoshenko beam]
The transmission line, which is modelled by the telegraph equations, provided with certain inputs and outputs formed by voltages and currents at the boundary (see \cite[Exercise~7.4]{Jacob2012}),
is also belonging to the class treated in this section. Our assumption on the Hamiltonian means that the longitudinal inductance $L:[a,b]\to\R$ as well as the transversal capacitance $C:[a,b]\to\R$ are positive-valued with $L,L^{-1},C,C^{-1}\in L^1([a,b];\R)$.\\
Likewise, for the Timoshenko beam (see \cite[Example~7.1.4]{Jacob2012}), the longitudinal mass density $\rho:[a,b]\to\R$, the rotary moment of a~cross section $I_\rho:[a,b]\to\R$, the product of Young's modulus and moment of interia $EI:[a,b]\to\R$, and the moment of interia of a~cross section $K:[a,b]\to\R$ can - in our setup - be assumed to be positive-valued with \[\rho,\rho^{-1},I_\rho,I_\rho^{-1},EI,(EI)^{-1},K,K^{-1}\in L^1([a,b];\R).\]
\end{rem}

\subsection{Maxwell's equations}

Let $\Omega\subset\R^3$ be a bounded Lipschitz domain with outward normal $\nu:\partial\Omega\to\R^3$, let $\epsilon,\mu\in L^\infty(\Omega;\C^{3\times 3})$ be pointwise Hermitian and positive definite with $\epsilon^{-1},\mu^{-1}\in L^\infty(\Omega;\C^{3\times 3})$, and let $g\in L^\infty(\Omega;\C^{3\times 3})$ such that $-g$ is pointwise dissipative.\\
We consider Maxwell's equations on $\Omega$, which - in compact form - read
\begin{subequations}\label{eq:Maxwell_syst}
\begin{equation}\label{eq:Maxwell}\begin{aligned}
\tfrac{\partial}{\partial t}\begin{pmatrix}\mathbf{B}(t,\xi)\\\mathbf{D}(t,\xi)\end{pmatrix}\!=&\begin{bmatrix}0&-\curl\\\curl&-g\end{bmatrix}\begin{pmatrix}\mu(\xi)^{-1}\mathbf{B}(t,\xi)\\\epsilon(\xi)^{-1}\mathbf{D}(t,\xi)\end{pmatrix},\quad%\;\begin{pmatrix}\mathbf{B}(0,\xi)\\\mathbf{D}(0,\xi)\end{pmatrix}\!=\!\begin{pmatrix}\mathbf{B}_0(\xi)\\\mathbf{D}_0(\xi)\end{pmatrix}&&
t\ge0, \xi\in\Omega\end{aligned}
\end{equation}
Hereby, $\mathbf{B}$ and $\mathbf{D}$ respectively stand for the magnetic flux and electric displacement, $\mu$ is the magnetic permeability, and $\epsilon$ is the electric permittivity. The quantities $\mathbf{H}:=\mu^{-1}\mathbf{B}$, $\mathbf{E}:=\epsilon^{-1}\mathbf{D}$ are the magnetic and electric field intensities, resp.\\
 Typically, Maxwell's equations are further provided with the conditions $\div \mathbf{H}=0$, $\div \mathbf{E}=\rho$ for some scalar field $\rho:\Omega\to\C$ (whose physical interpretation is that of charge density). By using a {\em Helmholtz decomposition}, that is, an orthogonal decomposition into a~divergence-free function and a~gradient field (and further using that the latter has trivial curl and is therefore constant in time), these conditions can be coded in the initial condition.\\
 As input, we choose a~tangential condition on the electric field, i.e.,
\begin{equation}\label{eq:maxinp}
u(t,\xi)=\nu(\xi)\times  \big(\epsilon(\xi)^{-1}\mathbf{D}(t,\xi)\big),\quad t\ge0, \xi\in\partial\Omega,
\end{equation}
where the outward normal fulfills $\nu\in L^\infty(\partial\Omega;\R^3)$ (with respect to the surface measure on $\partial\Omega$) as $\Omega$ is a~Lipschitz domain. Consider the function $\pi:\partial\Omega\to \R^{3\times 3}$ which maps $\xi\in\partial\Omega$ to the orthogonal projection onto the tangential space of $\partial\Omega$ at $\xi$. By using that
\[
\pi(\xi)w=(\nu(\xi)\times w)\times \nu(\xi)=w-(\nu(\xi)^\top w)\nu(\xi)\quad\forall\,w\in\C^3\text{ and almost all }\xi\in\partial\Omega,
\]
we have $\pi\in L^\infty(\partial\Omega;\R^{3\times 3})$. Our output is formed by the tangential component of the magnetic field. That is,
\[y(t,\xi)= {\pi(\xi)}\big(\mu(\xi)^{-1}\mathbf{B}(t,\xi)\big),\quad t\ge0,\, \xi\in\partial\Omega.\]
To show that this system fits into the framework of port-Hamiltonian system nodes, we next introduce the involved spaces. We take  $\calX=L^2(\Omega;\C^3)^2$, and the Hamiltonian is $\calH:\calX\to\R$ with
\begin{equation}\calH\spvek{\mathbf{B}}{\mathbf{D}}=\tfrac12\int_\Omega\mu(\xi)^{-1}\|\mathbf{B}(\xi)\|^2+\epsilon(\xi)^{-1}\|\mathbf{D}(\xi)\|^2{\rm d}\xi\label{eq:maxout}\end{equation}
\end{subequations}
that is, the sum of electric and magnetic energy. Our boundedness assumptions on the magnetic permeability and the electric permittivity yield that the sesquilinear form associated with $\calH$ is bounded and coercive. Hence, by Remark~\ref{rem:dissnode}\,\eqref{rem:dissnode4}, the energetic space and its anti-dual fulfill $\calX = \calX_h = \calX_h^* = (L^2(\Omega;\C^3))^2$. Similar as in the previous sections, these spaces are provided with the inner products
\begin{subequations}\label{eq:indualMax}
\begin{align}
\left\langle\spvek{\mathbf{B}_1}{\mathbf{D}_1},\spvek{\mathbf{B}_2}{\mathbf{D}_2}\right\rangle_{\calX_h}&=\int_\Omega \mathbf{B}_2(\xi)^*\mu(\xi)^{-1}\mathbf{B}_1(\xi)+\mathbf{D}_2(\xi)^*\epsilon(\xi)^{-1}\mathbf{D}_1(\xi){\rm d}\xi,\\
\left\langle\spvek{\mathbf{H}_1}{\mathbf{E}_1},\spvek{\mathbf{H}_2}{\mathbf{E}_2}\right\rangle_{\calX_h^*}&=\int_\Omega\mathbf{H}_2(\xi)^*\mu(\xi)\mathbf{H}_1(\xi)
+\mathbf{E}_2(\xi)^*\epsilon(\xi)\mathbf{E}_1(\xi){\rm d}\xi,
\end{align}
whereas the duality product is simply the inner product in $L^2(\Omega;\C^3)$, i.e.,
\begin{align}
\left\langle\spvek{\mathbf{H}}{\mathbf{E}},\spvek{\mathbf{B}}{\mathbf{D}}\right\rangle_{\calX_h^*,\calX_h}&=\int_\Omega \mathbf{B}(\xi)^*\mathbf{H}(\xi)+\mathbf{D}(\xi)^*\mathbf{E}(\xi){\rm d}\xi.
\end{align}
\end{subequations}
Note that all of the above three products have the physical dimension of energy.\\
The choice of suitable input and output spaces is more involved; this is subject of the following part which is inspired by \cite{Skrepek21,BuCoSh02,WeSt13}. We first introduce the space of tangential vector fields by
\[L^2_\tau(\partial\Omega):=\pi L^2(\partial\Omega;\C^3)=\setdef{w\in L^2(\partial\Omega;\C^3)}{\nu^\top w=0\text{ on }\partial\Omega},\]
and the {\em tangential component trace operator $\pi_\tau$}, which maps to the pointwise projection of the trace to the tangential space. In other words, we perform a~pointwise multiplication of the tangential projection with the boundary trace,
\[\begin{aligned}
\pi_\tau:\qquad W^{1,2}(\Omega;\C^3)&\to L^2_\tau(\partial\Omega),\\
\mathbf{A}&\mapsto \pi\,\big(\gamma\mathbf{A}\big),
\end{aligned}\]
where $\gamma:W^{1,2}(\Omega;\C^3)\to W^{1/2,2}(\partial\Omega;\C^3)$ is the trace operator. By the continuous embedding $W^{1/2,2}(\partial\Omega;\C^3)\subset L^2(\partial\Omega;\C^3)$ and $\pi\in L^\infty(\partial\Omega,\R^{3\times 3})$, we have \[\pi_\tau\in L(W^{1,2}(\Omega;\C^3),L^2_\tau(\partial\Omega)).\]
We now set
\[
V_\pi:=\pi_\tau W^{1,2}(\Omega;\C^3),
\]
which is equipped with the norm
\[
\|h\|_{V_\pi} := \inf\setdef{\|\mathbf{A}\|_{W^{1/2,2}(\Omega;\C^3)}}{\mathbf{A}\in W^{1,2}(\Omega)\,\wedge\,h=\pi_\tau \mathbf{A}}.
\]
The trace theorem \cite[Thm.~1.5.1.3]{Gris85} now yields that $\pi_\tau:W^{1,2}(\Omega)\to V_\pi$ is bounded, whereas surjectivity follows from defintion of $V_\pi$. Further, $V_\pi\subset L^2_\tau(\partial\Omega)$ and the tangential component trace operator can also be regarded as a~bounded mapping from $W^{1,2}(\Omega;\C^3)$ to $L^2_\tau(\partial\Omega)$.\\
The {\em tangential trace} of $\mathbf{A}\in W^{1,2}(\Omega;\C^3)$ is the pointwise cross product of the boundary trace of $\mathbf{A}$ with the outward normal $\nu\in L^\infty(\Omega;\R^3)$, i.e., $\widetilde{\gamma}_\tau\mathbf{A}:=\nu\times \gamma \mathbf{A}\in L^2_\tau(\partial\Omega)$, where the latter holds by $(\nu\times \gamma \mathbf{A})\,\bot\, \nu$ almost everywhere on $\partial\Omega$ (with respect to the surface measure on $\partial\Omega$). It has been shown in \cite[Prop.~4.3]{WeSt13} that $\widetilde{\gamma}_\tau$ uniquely extends to a~(not necessarily surjective) bounded mapping
\[\gamma_\tau: \qquad H(\curl,\Omega)\to V_\pi^*,\]
where $H(\curl,\Omega)$ consists of all elements of $L^2(\Omega;\C^3)$ whose weak curl is again in $L^2(\Omega;\C^3)$. This is a~Hilbert space equipped with the norm
\[\|\mathbf{A}\|_{H(\curl,\Omega)}=\big(\|\mathbf{A}\|_{L^2(\Omega;\C^3)}^2+\|\curl\mathbf{A}\|_{L^2(\Omega;\C^3)}^2\big)^{1/2}.\]
The basic ingredient for this is Green's formula
\begin{multline}\label{eq:curlgreen}
\forall \,\varphi\in W^{1,2}(\Omega;\C^3),\, \mathbf{A}\in H(\curl,\Omega):\\ \langle \curl\mathbf{A},\varphi\rangle_{L^2(\Omega;\C^3)}-\langle \mathbf{A},\curl\varphi\rangle_{L^2(\Omega;\C^3)}=\langle\gamma_\tau \mathbf{A},\pi_\tau \varphi\rangle_{V_\pi^*,V_\pi},
\end{multline}
where the latter expression restricts to an $L^2$-inner product $\langle\gamma_\tau \mathbf{A},\pi_\tau \varphi\rangle_{L^2(\partial\Omega;\C^3)}$,
%\[ \langle \curl\mathbf{E},\varphi\rangle_{L^2(\Omega;\C^3)}-\langle \mathbf{E},\curl\varphi\rangle_{L^2(\Omega;\C^3)}=\langle\gamma_\tau \mathbf{E},\pi_\tau \varphi\rangle_{L^2_\tau(\partial\Omega)},\]
if both $\varphi$ and $\mathbf{A}$ are in $W^{1,2}(\Omega;\C^3)$ \cite[Chap.~7]{BuCoSh02}. These spaces form indeed a~Gelfand triple with $V_\pi\subset L^2_\tau(\partial\Omega)\subset V_\pi^*$.\\
Now consider the inner product space
\begin{align*}
H(\curl,\Omega,\tau)
&=\setdef{\mathbf{A}\in H(\curl,\Omega)}{\,\gamma_\tau\mathbf{A}\in L^2_\tau(\partial\Omega)},\\
\|\mathbf{A}\|_{H(\curl,\Omega,\tau)}
&=\big(\|\mathbf{A}\|_{L^2(\Omega;\C^3)}^2 + \|\curl\mathbf{A}\|_{L^2(\Omega;\C^3)}^2 + \|\gamma_\tau\mathbf{A}\|_{L^2(\partial\Omega;\C^3)}^2\big)^{1/2},
\end{align*}
which is even a~Hilbert space as shown in \cite[Prop.~4.3]{WeSt13}. The construction of $H(\curl,\Omega,\tau)$ yields that the tangential trace operator can be regarded as a~bounded (even contractive) operator from $H(\curl,\Omega,\tau)$ to $L^2_\tau(\partial\Omega)$. Hence, we may consider the orthogonal projector $P_\tau\in L(H(\curl,\Omega,\tau))$ along $\ker\gamma_\tau$. We introduce the sesquilinear form
\begin{equation}\label{eq:htau}
\begin{aligned}
h_\tau:\qquad &\dom h_\tau\times \dom h_\tau\to\C,\qquad \dom h_\tau=\gamma_\tau H(\curl,\Omega,\tau),\\
&(f,g)\mapsto \langle P_\tau\mathbf{A}_f, P_\tau\mathbf{A}_g\rangle_{H(\curl,\Omega)}\quad \text{ for some $\mathbf{A}_f,\mathbf{A}_g\in H(\curl,\Omega,\tau)$}\\&\qquad\qquad\qquad\qquad\qquad\qquad\qquad\qquad\qquad\text{ with $f=
{\gamma_\tau}\mathbf{A}_f$, $g={\gamma_\tau}\mathbf{A}_g$}.
\end{aligned}
\end{equation}
It can be easily seen that this form is well-defined. The projection theorem \cite[Thm.~4.3]{Alt16} moreover implies that
\begin{multline}\label{eq:projmin}
\forall\, f\in \gamma_\tau H(\curl,\Omega,\tau):\\ h_\tau(f,f)=\inf\setdef{\|\mathbf{A}\|_{H(\curl,\Omega)}^2}{\mathbf{A}\in H(\curl,\Omega,\tau)\text{ with }\gamma_\tau\mathbf{A}=f}.
\end{multline}
To show that Maxwell equations with the above input-output configuration form a~port-Hamil\-tonian system in the sense of Section~\ref{sec:phsysnodes}, our next aim is to introduce the input and output spaces $\calU$, $\calU^*$.
Such spaces are provided in \cite{BuCoSh02} by showing that these are spaces of certain (fractional) regularity with respect to certain tangential differential operators. Our approach is somewhat less concrete, but also -- at least in our opinion -- more handy. We construct the input and output space from a~quasi Gelfand triple with respect to the pivot space $L^2_\tau(\partial\Omega)$, an approach similar to that in \cite{Skrepek21}. This will indeed been done on the basis of the sesquilinear form $h_\tau$.

\begin{prop}[Properties of $h_\tau$]
Let $\Omega\subset\R^3$ be a bounded Lipschitz domain. Then $h_\tau$ as in \eqref{eq:htau} is a densely defined and closed positive sesquilinear form on $L^2_\tau(\partial\Omega)$.
\end{prop}
\begin{proof}
Completeness of $H(\curl,\Omega,\tau)$ directly implies closedness of $h_\tau$, whereas positivity follows by construction of $h_\tau$. Density of $\dom h_\tau$ in $L^2_\tau(\partial\Omega)$ holds, since $W^{1,2}(\Omega;\C^3)$ is dense in $H(\curl,\Omega)$ and
\[
{\gamma}_\tau\mathbf{A} = \nu\times\gamma\mathbf{A}\in L^2_\tau(\partial\Omega)\quad \forall\,\mathbf{A}\in W^{1,2}(\Omega;\C^3).
\]
\end{proof}

We define the input space $\calU$ to be the energetic space induced by the form $h_\tau$, cf.\ Definition \ref{d:energetic}. Hence, $(\calU^*,L^2_\tau(\partial\Omega),\calU)$ is the quasi Gelfand triple associated with $h_\tau$. By using $W^{1,2}(\Omega;\C^3)\subset H(\curl,\Omega,\tau)\subset H(\curl,\Omega)$, we obtain that $H(\curl,\Omega,\tau)$ is dense in $H(\curl,\Omega)$. Further, \eqref{eq:projmin} gives rise to
\[
\forall\mathbf{A}\in H(\curl,\Omega,\tau):\quad\|\gamma_\tau\mathbf{A}\|^2_{\calU}=h_\tau(\gamma_\tau\mathbf{A},\gamma_\tau\mathbf{A})\leq \|\mathbf{A}\|^2_{H(\curl,\Omega)}.
\]
Hence, the tangential trace operator extends to a~bounded operator from $H(\curl,\Omega)$ to $\calU$, which is surjective by construction of $\calU$, and in the following also denoted by $\gamma_\tau$. Hence, the input equation \eqref{eq:maxinp} will be rewritten by
\[
\gamma_\tau\big(\epsilon^{-1}\mathbf{D}(t)\big)=u(t)\in\calU,\qquad t\geq0.
\]
To properly rewrite the output equation, we are motivated by \eqref{eq:curlgreen} to redefine the tangential component trace operator by $\pi_\tau:H(\curl,\Omega)\to\calU^*$ with
\begin{equation}
\label{eq:curlgreen_gen}
\forall \,\mathbf{A},\mathbf{C}\in H(\curl,\Omega):\\
\langle \pi_\tau\mathbf{A},\gamma_\tau \mathbf{C}\rangle_{\calU^*,\calU}
:=\langle \mathbf{A},\curl\mathbf{C}\rangle_{L^2(\Omega;\C^3)}-\langle \curl\mathbf{A},\mathbf{C}\rangle_{L^2(\Omega;\C^3)}.
\end{equation}
Surjectivity of $\gamma_\tau:H(\curl,\Omega)\to\calU$ yields that $\pi_\tau$ is well-defined. Moreover, by \eqref{eq:curlgreen},
$\gamma_\tau\mathbf{A}$ is truly the ``well-tried tangential component trace'', if $\mathbf{A}\in W^{1,2}(\Omega;\C^3)$.
Based on this, the output equation \eqref{eq:maxout} will be rewritten by
\[
\pi_\tau\big(\mu^{-1}\mathbf{B}(t)\big)=y(t)\in\calU^*,\qquad t\geq0.
\]
Now we are able to introduce the dissipation node that corresponds to the system \eqref{eq:Maxwell_syst}.

\begin{prop}[Dissipation node for Maxwell's equations]\label{prop:maxwdiss}
Let $\Omega\subset\R^3$ be a~Lipschitz domain, let $\epsilon,\mu\in L^\infty(\Omega;\C^{3\times 3})$ be pointwise Hermitian and positive definite with $\epsilon^{-1},\mu^{-1}\in L^\infty(\Omega;\C^{3\times 3})$, and let $g\in L^\infty(\Omega;\C^{3\times 3})$ such that $-g$ is pointwise dissipative.
Further, let $(\calU^*,L^2_\tau(\partial\Omega),\calU)$ be the quasi Gelfand triple associated with the sesquilinear form $h_\tau$ as in \eqref{eq:htau}. Then, for $\calX_h = \calX_h^* = (L^2(\Omega;\C^3))^2$ with inner products as in \eqref{eq:indualMax}, the operator $ M:\calX_h^*\times\calU\supset\dom M\to\calX_h\times\calU^*$ with
\[
\dom M = \setdef{\left(\begin{smallmatrix}\mathbf{H}\\\mathbf{E}\\u\end{smallmatrix}\right)\in H(\curl,\Omega)^2\times \calU}{\,u=\gamma_\tau\mathbf{E}},\qquad
 M\left(\begin{smallmatrix}\mathbf{H}\\\mathbf{E}\\u\end{smallmatrix}\right)=\left[\begin{smallmatrix}0&-\curl\\\curl&-g\\-\pi_\tau &0\end{smallmatrix}\right]\spvek{\mathbf{H}}{\mathbf{E}}
\]
is a~dissipation node on $(\calX_h,\calU)$.
\end{prop}
\begin{proof}
We successively verify that $M$ fulfills the criteria \eqref{def:dissnode1}--\eqref{def:dissnode4} in Definition~\ref{def:dissnode}. By Remark~\ref{rem:dissnode}\,\eqref{rem:dissnode4}, it is no loss of generality to assume that $\calX_h$ and $\calX_h^*$ are both provided with the standard $L^2$-inner product (equivalently, $\mu=\epsilon\equiv1$).
\begin{enumerate}[(a)]
\item\label{prop:maxwdiss1} We start with showing that ${M}$ is dissipative. By using the definition of the tangential component trace operator in \eqref{eq:curlgreen_gen}, we have
\begin{align*}
\forall\,\left(\begin{smallmatrix}\mathbf{H}\\\mathbf{E}\\u\end{smallmatrix}\right)\in \dom M:
&\phantom{=}\Re\left\langle\left(\begin{smallmatrix}\mathbf{H}\\\mathbf{E}\\u\end{smallmatrix}\right), M\left(\begin{smallmatrix}\mathbf{H}\\\mathbf{E}\\u\end{smallmatrix}\right)\right\rangle_{\calX\times\calU,\calX\times\calU^*}\\
&{=}\Re\left\langle\left(\begin{smallmatrix}\mathbf{H}\\\mathbf{E}\\\gamma_\tau\mathbf{E}\end{smallmatrix}\right),\left(\begin{smallmatrix}-\curl\mathbf{E}\\\curl\mathbf{H}-g\mathbf{E}\\-\pi_\tau \mathbf{H}\end{smallmatrix}\right)\right\rangle_{\calX\times\calU,\calX\times\calU^*}\\
%&{=}-\Re\langle \mathbf{E},g\mathbf{E}\rangle_{L^2(\Omega;\C^3)}
%-\Re\langle \mathbf{H},\curl\mathbf{E}\rangle_{L^2(\Omega;\C^3)}+\Re\langle \mathbf{E},\curl\mathbf{H}\rangle_{L^2(\Omega;\C^3)}
%-\Re\langle \pi_\tau\mathbf{E},\gamma_\tau \mathbf{H}\rangle_{\calU^*,\calU}\\
&{=}-\Re\langle \mathbf{E},g\mathbf{E}\rangle_{L^2(\Omega;\C^3)}\\&\quad
+\Re\left(\langle \mathbf{E},\curl\mathbf{H}\rangle_{L^2(\Omega;\C^3)}-\langle \curl\mathbf{E},\mathbf{H}\rangle_{L^2(\Omega;\C^3)}-
\langle \pi_\tau\mathbf{E},\gamma_\tau \mathbf{H}\rangle_{\calU^*,\calU}\right)\\
&{=}-\Re\langle \mathbf{E},g\mathbf{E}\rangle_{L^2(\Omega;\C^3)}\leq0.
\end{align*}
Next we show that $ M$ is closed. Assume that
\[\left(\begin{smallmatrix}\mathbf{H}_n\\\mathbf{E}_n\\u_n\end{smallmatrix}\right)\to \left(\begin{smallmatrix}\mathbf{H}\\\mathbf{E}\\u\end{smallmatrix}\right)\text{ in }L^2(\Omega;\C^3)^2\times\calU,\qquad M\left(\begin{smallmatrix}\mathbf{H}_n\\\mathbf{E}_n\\u_n\end{smallmatrix}\right)\to \left(\begin{smallmatrix}\mathbf{A}\\\mathbf{C}\\w\end{smallmatrix}\right)\text{ in }L^2(\Omega;\C^3)^2\times\calU^{*}.\]
By definition of $M$, we obtain that both $(\mathbf{H}_n)$ and $(\mathbf{E}_n)$ are Cauchy sequences in $H(\curl,\Omega)$, and thus convergent. Then, by $L^2$-convergence towards $\mathbf{H}$ and $\mathbf{E}$, resp., the sequences $(\mathbf{H}_n)$ and $(\mathbf{E}_n)$ converge in $H(\curl,\Omega)$ to $\mathbf{H}$ and $\mathbf{E}$, resp. This in particular shows that $(\curl\mathbf{H}_n)$ and $(\curl\mathbf{E}_n)$ converge in $L^2(\Omega;\C^3)$ to $\curl\mathbf{H}$ and $\curl\mathbf{E}$, resp. A~direct consequence is that $\mathbf{A}=-\curl\mathbf{E}$ and $\mathbf{C}=\curl\mathbf{H}-g\mathbf{E}$.
Since, moreover, boundedness of $\gamma_\tau:H(\curl,\Omega)\to\calU$ and $\pi_\tau:H(\curl,\Omega)\to\calU^*$ give
\begin{align*}
u\stackrel{n\to\infty}{\longleftarrow}&u_n=\gamma_\tau \mathbf{E}_n\stackrel{n\to\infty}{\longrightarrow}\gamma_\tau \mathbf{E},\\
w\stackrel{n\to\infty}{\longleftarrow}&-\pi_\tau \mathbf{H}_n\stackrel{n\to\infty}{\longrightarrow}-\pi_\tau \mathbf{H},
\end{align*}
we are led to
\[\left(\begin{smallmatrix}\mathbf{H}\\\mathbf{E}\\u\end{smallmatrix}\right)\in\dom M,\qquad M\left(\begin{smallmatrix}\mathbf{H}\\\mathbf{E}\\u\end{smallmatrix}\right)=\left(\begin{smallmatrix}\mathbf{A}\\\mathbf{C}\\w\end{smallmatrix}\right).\]
\item\label{prop:maxwdiss2} This can be achieved by exactly the same argumentation as in \eqref{prop:maxwdiss1}.
\item\label{prop:maxwdiss3} Let $u\in\calU$. By surjectivity of the tangential trace operator $\gamma_\tau:H(\curl;\Omega)\to\calU$, there exists some $\mathbf{E}\in H(\curl,\Omega)$ with $u=\gamma_\tau\mathbf{E}$. Then
\[\left(\begin{smallmatrix}0\\\mathbf{E}\\u\end{smallmatrix}\right)\in\dom M.\]
\item\label{prop:maxwdiss4} Consider the main operator
\[F:\dom F=\setdef{\spvek{\mathbf{H}}{\mathbf{E}}\in H(\curl;\Omega)^2}{\gamma_\tau \mathbf{E}=0},\quad F=\left[\begin{smallmatrix}0&-\curl\\\curl&-g\end{smallmatrix}\right].
\]
{\em Step~1:} We show that $\dom F^*=\dom F$ with
\[F^*=\left[\begin{smallmatrix}0&\curl\\-\curl&-g^*\end{smallmatrix}\right].\]
%where $g^*$ stands for the pointwise conjugate transpose of $g$.\\
For $\spvek{\mathbf{H}}{\mathbf{E}},\spvek{\mathbf{A}}{\mathbf{C}}\in\dom F$, the definition of the tangential component trace operator in \eqref{eq:curlgreen_gen} yields
\begin{align*}
&\phantom{=}\left\langle \spvek{\mathbf{A}}{\mathbf{C}},F\spvek{\mathbf{H}}{\mathbf{E}}\right\rangle_{L^2(\Omega;\C^3)^2}\\
&{=}\left\langle \spvek{\mathbf{A}}{\mathbf{C}},\spvek{-\curl\mathbf{E}}{\curl\mathbf{H}-g\mathbf{E}}\right\rangle_{L^2(\Omega;\C^3)^2}\\
&=-\langle\mathbf{A},\curl\mathbf{E}\rangle_{L^2(\Omega;\C^3)}+\langle\mathbf{C},\curl\mathbf{H}\rangle_{L^2(\Omega;\C^3)}-\langle\mathbf{C},g\mathbf{E}\rangle_{L^2(\Omega;\C^3)}\\
&=-\langle\curl\mathbf{A},\mathbf{E}\rangle_{L^2(\Omega;\C^3)}+\langle\curl\mathbf{C},\mathbf{H}\rangle_{L^2(\Omega;\C^3)}-\langle g^*\mathbf{C},\mathbf{E}\rangle_{L^2(\Omega;\C^3)}\\
&{=}\left\langle \spvek{\curl\mathbf{C}}{-\curl\mathbf{A}-g^*\mathbf{C}},\spvek{\mathbf{H}}{\mathbf{E}}\right\rangle_{L^2(\Omega;\C^3)^2}.
\end{align*}
It remains to show that $\dom F^*\subset\dom F$. Assume that $\spvek{\mathbf{A}}{\mathbf{C}}\in\dom F^*$. Then there exists some $\spvek{\mathbf{F}}{\mathbf{G}}\in L^2(\Omega;\C^3)$, such that for all $\spvek{\mathbf{H}}{\mathbf{E}}\in \dom F$,
\begin{align*}
\left\langle \spvek{\mathbf{F}}{\mathbf{G}},\spvek{\mathbf{H}}{\mathbf{E}}\right\rangle_{L^2(\Omega;\C^3)^2}=\left\langle \spvek{\mathbf{A}}{\mathbf{C}},F\spvek{\mathbf{H}}{\mathbf{E}}\right\rangle_{L^2(\Omega;\C^3)^2}{=}\left\langle \spvek{\mathbf{A}}{\mathbf{C}},\spvek{-\curl\mathbf{E}}{\curl\mathbf{H}-g\mathbf{E}}\right\rangle_{L^2(\Omega;\C^3)^2}.
\end{align*}
By first setting $\mathbf{E}=0$, we obtain
\begin{align*}
\forall \,\mathbf{H}\in H(\curl,\Omega):\quad \langle\mathbf{F},\mathbf{H}\rangle_{L^2(\Omega;\C^3)}=\langle\mathbf{C},\curl\mathbf{H}\rangle_{L^2(\Omega;\C^3).}
\end{align*}
By first choosing compact supported smooth functions for $\mathbf{H}$, the definition of the weak curl gives rise to  $\mathbf{F}=\curl\mathbf{C}$, and the definition of the tangential component trace operator in \eqref{eq:curlgreen_gen} yields
\begin{multline*}
\forall \,\mathbf{H}\in H(\curl,\Omega):\\ \langle\curl\mathbf{C},\mathbf{H}\rangle_{L^2(\Omega;\C^3)}=\langle\mathbf{C},\curl\mathbf{H}\rangle_{L^2(\Omega;\C^3)}=\langle\curl\mathbf{C},\mathbf{H}\rangle_{L^2(\Omega;\C^3).}
+\langle \pi_\tau\mathbf{C},\gamma_\tau \mathbf{H}\rangle_{\calU^*,\calU}.
\end{multline*}
This leads to $\pi_\tau\mathbf{C}=0$, which implies, as shown in \cite{BuCoSh02}, that $\gamma_\tau\mathbf{C}=0$.\\
Now setting $\mathbf{H}=0$, we obtain
\[
\forall \,\mathbf{E}\in H(\curl,\Omega)\text{ with }\gamma_\tau\mathbf{E}=0: \langle\mathbf{G}+g^*\mathbf{C},\mathbf{E}\rangle_{L^2(\Omega;\C^3)}=-\langle\mathbf{A},\curl\mathbf{E}\rangle_{L^2(\Omega;\C^3)},
\]
by again choosing compactly supported smooth functions for $\mathbf{E}$, the definition of the weak curl gives rise to  $\mathbf{G}+g^*\mathbf{C}=-\curl\mathbf{A}$. This altogether shows that $\spvek{\mathbf{A}}{\mathbf{C}}\in\dom F$.\\
{\em Step~2:} We show that $F$ has property {\eqref{def:dissnode4}} in Definition~\ref{def:dissnode}.\\
By using step~1, we have for all $\spvek{\mathbf{A}}{\mathbf{C}}\in\dom F^*=\dom F$
\begin{align*}
\left\langle \spvek{\mathbf{A}}{\mathbf{C}},F^*\spvek{\mathbf{A}}{\mathbf{C}}\right\rangle_{L^2(\Omega;\C^3)^2}
&{=}\left\langle \spvek{\mathbf{A}}{\mathbf{C}},\spvek{\curl\mathbf{C}}{-\curl\mathbf{A}-g^*\mathbf{C}}\right\rangle_{L^2(\Omega;\C^3)^2}\\
&{=}\langle \mathbf{A},\curl\mathbf{C}\rangle_{L^2(\Omega;\C^3)}-\langle \mathbf{C},\curl\mathbf{A}\rangle_{L^2(\Omega;\C^3)}-\langle \mathbf{C},g^*\mathbf{C}\rangle_{L^2(\Omega;\C^3)}\\
&{=}-\langle \mathbf{C},g^*\mathbf{C}\rangle_{L^2(\Omega;\C^3)}\leq0.
\end{align*}
Hence, $F^*$ is dissipative, and we can conclude the desired result from Proposition~\ref{prop:maxdiss}.
\end{enumerate}
\end{proof}

\section{Summary}

We have introduced an~operator theoretic approach to linear partial differential equations with port-Hamiltonian structure, that in particular covers the case of boundary control and observation. Our basic ingredient is the theory of ``system nodes'' by {\sc Staffans} in \cite{Staffans2005} and antecedent journal publications by the same author. Our presented approach also goes along with a~solution theory, and it allows to show that the solutions fulfill a~dissipation inequality. The latter has -- in physical terms -- the meaning of an energy balance. The theory is applied to an advection-diffusion equation and Maxwell's equations, both on bounded Lipschitz domains and subject to boundary control and observation. Further, the class of hyperbolic systems on one-dimensional spatial domains as treated in the book \cite{Jacob2012} by {\sc Jacob} and {\sc Zwart} has been analyzed in the context of the presented theory, and further generalized to systems with certain unbounded and/or non-coercive Hamiltonian densities.\\
The authors believe that a~variety of further physically motivated partial differential equations fit into the presented framework. Besides taking further application fields into account (such as, for instance, from fluid dynamics), there are rather many open research questions that go beyond the content of this article. For instance, such topics include interconnection of port-Hamiltonian systems as well as relating the presented theory to the one in \cite{van2021differential,van2002hamiltonian,van2014port} on (Stokes-)Dirac and resistive structures.

\section*{Acknowledgement}

%While their research on this topic was at a~``green banana stage'', the authors gave a~couple of talks on the state of affairs, such that it could ripen at the customers' places.
The authors would like to express their gratitude to Nathanael Skrepek (TU BA Freiberg) for his very fruitful comments on the choice of boundary trace spaces for Maxwell's equations and for letting them know about the concept of quasi Gelfand triples. Also thanks to Hans Zwart (U~Twente) for providing references and valuable insight into the history of the systems treated in Section~\ref{sec:bjhz}.

\bibliographystyle{abbrv}

\end{document}